\documentclass{amsart}
\usepackage{latexsym}
\usepackage{amsmath, amssymb}
\usepackage[all]{xy}
\newtheorem{dfs}{Definition}[section]
\newtheorem{lms}[dfs]{Lemma}
\newtheorem{thms}[dfs]{Theorem}
\newtheorem{props}[dfs]{Proposition}
\newtheorem{qus}[dfs]{Question}
\newtheorem{cors}[dfs]{Corollary}
\newtheorem{rems}[dfs]{Remark}

\begin{document}

\title[Stability in the Cuntz semigroup]{Stability in the Cuntz semigroup of a commutative $C^*$-algebra}

\author{Andrew S. Toms}

\maketitle

\begin{abstract}
Let $A$ be a $C^*$-algebra.
The Cuntz semigroup $W(A)$ is an analogue for positive elements of the semigroup 
$V(A)$ of Murray-von Neumann equivalence classes of projections in matrices over $A$.
We prove stability theorems for the Cuntz semigroup of a commutative $C^*$-algebra
which are analogues of classical stability theorems for topological vector bundles over 
compact Hausdorff spaces.  

Let $\mathcal{SDG}$ denote the class of simple, unital, and infinite-dimensional
AH algebras with slow dimension growth, and let $A$ be an element of $\mathcal{SDG}$.  We apply our stability
theorems to obtain the following:

\begin{enumerate}
\item[(i)]   $A$ has strict comparison of positive elements;
\item[(ii)]  $W(A)$ is recovered functorially from the Elliott invariant of $A$;
\item[(iii)]  The lower semicontinuous dimension functions on $A$ are weak-$*$ dense 
in the dimension functions on $A$;
\item[(iv)]  The dimension functions on $A$ form a Choquet simplex.
\end{enumerate}
Statement (ii) confirms a conjecture of Perera and the author, while statements (iii) and
(iv) confirm, for $\mathcal{SDG}$, conjectures of Blackadar and Handelman from the early 1980s.
\end{abstract}

\section{Introduction}\label{intro}
In 1978, Cuntz introduced a generalisation of Murray-von Neumann comparison:
given positive elements $a$ and $b$ in a $C^*$-algebra $A$, write $a \precsim b$ if
there is a sequence $(v_i)_{i=1}^{\infty}$ in $A$ such that
\[
||v_i b v_i^* - a || \stackrel{i \to \infty}{\longrightarrow} 0
\]
(\cite{Cu1}).  We say that $a$ is {\it Cuntz subequivalent} to $b$.
The relation $\sim$ given by
\[
a \sim b \Leftrightarrow a \precsim b \ \wedge \ b \precsim a
\]
is an equivalence relation known as {\it Cuntz equivalence}.  If $A$ is unital, 
then one can mimic the construction of the ordered semigroup 
$V(A)$ of Murray-von Neumann equivalence classes of projections in matrices over $A$
by substituting positive elements for projections and Cuntz equivalence
for Murray-von Neumann equivalence.  This yields a positively ordered Abelian monoid
$W(A)$ called the {\it Cuntz semigroup} of $A$, and its partially ordered Grothendieck
envelope $\mathrm{K}_0^*(A)$.  If $A$ is stably finite, then each tracial state on
$A$ gives rise to an order preserving state on $W(A)$.  If the partial order
on $W(A)$ is determined by these states, then $A$ is said to have {\it strict comparison of
positive elements}, or simply {\it strict comparison}.

The study of the Cuntz semigroup has been resurgent of late.
R{\o}rdam has proved that simple, unital, exact, and finite 
$C^*$-algebras which absorb the Jiang-Su algebra tensorially have
strict comparsion (\cite{R4}), whence, by results of W. Winter and the author,
this last property is enjoyed by all classes of simple, separable, and
nuclear $C^*$-algebras currently known to satisfy Elliott's classification
conjecture (\cite{TW1},\cite{TW2}).  Coward, Elliott, and Ivanescu have 
recently identified a category of partially ordered semigroups into 
which the Cuntz semigroup is a continuous functor with respect to inductive limits.
Significantly, the Cuntz semigroup is deeply connected to the classification programme for separable and
nuclear $C^*$-algebras:  such algebras cannot be classified 
by their $\mathrm{K}$-theory in the absence of strict comparison (\cite{To2}), 
and there is some evidence that the converse will hold (\cite{TW3}).  Brown, Perera,
and the author recently proved a structure theorem for the Cuntz semigroup
which applies to most of our stock-in-trade simple, separable, and nuclear $C^*$-algebras
(\cite{BPT}, \cite{PT}).

In this paper we study the structure of the Cuntz semigroup for 
commutative and approximately homogeneous (AH) $C^*$-algebras.  Our main result is a positive answer 
to the following question, raised in \cite{To3}:
\begin{qus}\label{mainq}
Does there exist a constant $K > 0$ such that for any compact metrisable Hausdorff space $X$, 
$n \in \mathbb{N}$, and positive elements $a,b \in \mathrm{M}_{n}(\mathrm{C}(X))$ for which
\[
\mathrm{rank}(a)(x) + K \mathrm{dim}(X) \leq \mathrm{rank}(b)(x), \ \forall x \in X,
\]
one has $a \precsim b$? 
\end{qus}
\noindent
This question asks for an analogue in positive elements of the following well known stability
theorem for vector bundles.
\begin{thms}[Theorems 1.2 and 1.5, Chapter 8, \cite{H}]\label{kstab}
Let $X$ be a compact metrisable Hausdorff space, and let $\omega$ and $\xi$ be complex vector bundles
over $X$.  If the fibre dimension of $\omega$ exceeds the fibre dimension of $\xi$ by 
at least $\lceil \mathrm{dim}(X)/2 \rceil$ over every point in $X$, then $\xi$ is isomorphic
to a sub-bundle of $\omega$.
\end{thms}
\noindent

Our positive answer to Question \ref{mainq} has several applications to AH algebras.
Recall that an AH algebra is an inductive limit $C^*$-algebra
$A= \lim_{i \to \infty}(A_i,\phi_i)$, where
\begin{equation}\label{decomp}
A_i = \bigoplus_{l=1}^{n_i} p_{i,l}(\mathrm{C}(X_{i,l}) \otimes \mathcal{K})p_{i,l}
\end{equation}
for compact connected Hausdorff spaces $X_{i,l}$, projections $p_{i,l} \in
\mathrm{C}(X_{i,l}) \otimes \mathcal{K}$, and natural numbers $n_i$.  If
$A$ is separable, then one may assume that the $X_{i,l}$ are finite CW-complexes
(\cite{Bl2}, \cite{Go1}).  The algebras $A_i$ are called {\it semi-homogeneous}, and the
inductive system $(A_i,\phi_i)$ is referred to as a {\it decomposition} for $A$.
All AH algebras in this paper are assumed to be separable.

If an AH algebra $A$ admits a decompostion as in (\ref{decomp}) for which 
\[
\mathrm{max}_{1 \leq l \leq n_i} \ \left\{ \frac{\mathrm{dim}(X_{i,1})}{\mathrm{rank}(p_{i,1})}, \ldots,
\frac{\mathrm{dim}(X_{i,n_i})}{\mathrm{rank}(p_{i,n_i})} \right\} \stackrel{i \to \infty}{\longrightarrow} 0,
\]
then we say that $A$ has {\it slow dimension growth}.  Let $\mathcal{SDG}$ denote the
class of simple, unital, and infinite-dimensional AH algebras with slow dimension
growth.  $\mathcal{SDG}$ was studied intensively during the early- to mid-1990s by Blackadar,
Bratteli, D\u{a}d\u{a}rlat, Elliott, Gong, Kumjian, Li, R{\o}rdam, Thomsen, and others
(see \cite{BBEK}, \cite{BDR}, \cite{BE}, \cite{DNNP}, \cite{D}, \cite{EG}, \cite{EGL2}, \cite{EGL}, and \cite{G2}). 
The crowning achievement of this study was the confirmation of Elliott's classification
conjecture for $\mathcal{SDG}$ under the additional assumption of real rank zero (\cite{D}, \cite{EG}, \cite{G2}),
and the same confirmation in the real rank one case under the stronger hypothesis
of {\it very slow dimension growth} for $A$ (\cite{EGL}):
\[
\mathrm{max}_{1 \leq l \leq n_i} \ \left\{ \frac{\mathrm{dim}(X_{i,1})^3}{\mathrm{rank}(p_{i,1})}, \ldots,
\frac{\mathrm{dim}(X_{i,n_i})^3}{\mathrm{rank}(p_{i,n_i})} \right\} \stackrel{i \to \infty}{\longrightarrow} 0.
\]
This strengthened hypothesis is thought by experts to be unneccessary, but there
has been no progress on this problem.  In fact, there has been no progress on the basic
structure of elements of $\mathcal{SDG}$ since the late 1990s.

In this paper we obtain significant new results on the structure slow dimension
growth AH algebras.  We use our positive answer to Question \ref{mainq} to prove 
that if $A \in \mathcal{SDG}$, then $A$ has strict comparison of positive elements.  There
is evidence that this powerful $\mathrm{K}$-theoretic condition will characterise those
simple, separable, and nuclear $C^*$-algebras which are amenable to $\mathrm{K}$-theoretic classification ---
it already does so for a class of AH algebras intersecting substantially with $\mathcal{SDG}$ (\cite{TW3}) --- and
so our result supports the belief that $\mathcal{SDG}$ will satisfy the Elliott conjecture
without the assumption of very slow dimension growth.  
By appealing to our recent work
with Brown and Perera (\cite{BPT}), we confirm several conjectures pertaining to the 
structure of the Cuntz semigroup for $A \in \mathcal{SDG}$.  First, $W(A)$ is recovered
functorially from the Elliott invariant of $A$, allowing one to attack the Elliott conjecture
for $\mathcal{SDG}$ with considerably more information than would appear to be contained in
the Elliott invariant alone.  (This functorial recovery was conjectured for a class of algebras
containing $\mathcal{SDG}$ by Perera and the author in \cite{PT}.)  Second, the states on $W(A)$
coming from traces on $A$ are weak-$*$ dense in the state space of $W(A)$.  Third, the state
space of $W(A)$ is a Choquet simplex.  The results concerning the state space of $W(A)$ were
conjectured by Blackadar and Handelman in 1982 (\cite{BH}), but were only known to hold for
commutative $C^*$-algebras at the time. 

In a separate paper we will apply our positive answer to Question \ref{mainq} to
answer, affirmatively, a question of N. C. Phillips:  ``Do there exist simple, separable,
nuclear, and non-$\mathcal{Z}$-stable C$^*$-algebras which agree on the Elliott invariant,
yet are not isomorphic?'' (See \cite{To4}.)

Our paper is orgainsed as follows:  Section \ref{prelim} contains the definition 
of the Cuntz semigroup and recalls some essential facts about Cuntz subequivalence; 
Section \ref{stab} contains the stability theorems which answer Question \ref{mainq};
in Section \ref{ah} we establish results (i)-(iv) of the abstract.

\vspace{3mm}
\noindent
{\it Acknowledgements.}  We thank N. Christopher Phillips and Wilhelm Winter for
several inspiring conversations had during a visit to M{\"u}nster in August of 2005
which incited us to write this paper.  We also thank Mikael R{\o}rdam for supplying
Lemma \ref{posapprox1}, and Joachim Cuntz for supporting our visit to M{\"u}nster.

\section{Preliminaries}\label{prelim}

Let $A$ be a $C^*$-algebra, and let $\mathrm{M}_n(A)$ denote the $n \times n$ 
matrices whose entries are elements of $A$.  If $A = \mathbb{C}$, then we simply write $\mathrm{M}_n$.
Let $\mathrm{M}_{\infty}(A)$ denote the algebraic limit of the
direct system $(\mathrm{M}_n(A),\phi_n)$, where $\phi_n:\mathrm{M}_n(A) \to \mathrm{M}_{n+1}(A)$
is given by
\[
a \mapsto \left( \begin{array}{cc} a & 0 \\ 0 & 0 \end{array} \right).
\]
Let $\mathrm{M}_{\infty}(A)_+$ (resp. $\mathrm{M}_n(A)_+$)
denote the positive elements in $\mathrm{M}_{\infty}(A)$ (resp. $\mathrm{M}_n(A)$). 
Given $a,b \in \mathrm{M}_{\infty}(A)_+$, we say that $a$ is {\it Cuntz subequivalent} to
$b$ (written $a \precsim b$) if there is a sequence $(v_n)_{n=1}^{\infty}$ of
elements of $\mathrm{M}_{\infty}(A)$ such that
\[
||v_nbv_n^*-a|| \stackrel{n \to \infty}{\longrightarrow} 0.
\]
We say that $a$ and $b$ are {\it Cuntz equivalent} (written $a \sim b$) if
$a \precsim b$ and $b \precsim a$.  This relation is an equivalence relation,
and we write $\langle a \rangle$ for the equivalence class of $a$.  The set
\[
W(A) := \mathrm{M}_{\infty}(A)_+/ \sim
\] 
becomes a positively ordered Abelian monoid when equipped with the operation
\[
\langle a \rangle + \langle b \rangle = \langle a \oplus b \rangle
\]
and the partial order
\[
\langle a \rangle \leq \langle b \rangle \Leftrightarrow a \precsim b.
\]
In the sequel, we refer to this object as the {\it Cuntz semigroup} of $A$.  The
Grothendieck enveloping group of $W(A)$ is denoted $\mathrm{K}_0^*(A)$.

Given $a \in \mathrm{M}_{\infty}(A)_+$ and $\epsilon > 0$, we denote by 
$(a-\epsilon)_+$ the element of $C^*(a)$ corresponding (via the functional
calculus) to the function
\[
f(t) = \mathrm{max}\{0,t-\epsilon\}, \ t \in \sigma(a).
\]
(Here $\sigma(a)$ denotes the spectrum of $a$.)  
The proposition below collects some facts about Cuntz subequivalence due 
to Kirchberg and R{\o}rdam.

\begin{props}[Kirchberg-R{\o}rdam (\cite{KR}), R{\o}rdam (\cite{R4})]\label{basics}
Let $A$ be a $C^*$-algebra, and $a,b \in A_+$.  
\begin{enumerate}
\item[(i)] $(a-\epsilon)_+ \precsim a$ for every $\epsilon > 0$.
\item[(ii)] The following are equivalent:
\begin{enumerate}
\item[(a)] $a \precsim b$;
\item[(b)] for all $\epsilon > 0$, $(a-\epsilon)_+ \precsim b$;
\item[(c)] for all $\epsilon > 0$, there exists $\delta > 0$ such that $(a-\epsilon)_+ \precsim (b-\delta)_+$.
\end{enumerate}
\item[(iii)] If $\epsilon>0$ and $||a-b||<\epsilon$, then $(a-\epsilon)_+ \precsim b$.
\end{enumerate}
\end{props} 

Now suppose that $A$ is unital and stably finite, and denote by $\mathrm{QT}(A)$
the space of normalised 2-quasitraces on $A$ (v. \cite[Definition II.1.1]{BH}).
Let $S(W(A))$ denote the set of additive and order preserving maps $s$ from $W(A)$ to $\mathbb{R}^+$
having the property that $s(\langle 1_A \rangle) = 1$.
Such maps are called {\it states}.  Given $\tau \in \mathrm{QT}(A)$, one may 
define a map $s_{\tau}:\mathrm{M}_{\infty}(A)_+ \to \mathbb{R}^+$ by
\begin{equation}\label{ldf}
s_{\tau}(a) = \lim_{n \to \infty} \tau(a^{1/n}).
\end{equation}
This map is lower semicontinous, and depends only on the Cuntz equivalence class
of $a$.  It moreover has the following properties:
\vspace{2mm}
\begin{enumerate}
\item[(i)] if $a \precsim b$, then $s_{\tau}(a) \leq s_{\tau}(b)$;
\item[(ii)] if $a$ and $b$ are mutually orthogonal, then $s_{\tau}(a+b) = s_{\tau}(a)+s_{\tau}(b)$;
\item[(iii)] $s_{\tau}((a-\epsilon)_+) \nearrow s_{\tau}(a)$ as $\epsilon \to 0$. 
\end{enumerate}
\vspace{2mm}
Thus, $s_{\tau}$ defines a state on $W(A)$.
Such states are called {\it lower semicontinuous dimension functions}, and the set of them 
is denoted $\mathrm{LDF}(A)$. $\mathrm{QT}(A)$ is a simplex (\cite[Theorem II.4.4]{BH}), 
and the map from $\mathrm{QT}(A)$ to $\mathrm{LDF}(A)$ defined by (\ref{ldf}) is 
bijective and affine (\cite[Theorem II.2.2]{BH}).  A {\it dimension function} on
$A$ is a state on $\mathrm{K}_0^*(A)$, assuming that the latter has been equipped
with the usual order coming from the Grothendieck map.  The set of dimension functions
is denoted $\mathrm{DF}(A)$.  $\mathrm{LDF}(A)$ is a (generally proper) face of
$\mathrm{DF}(A)$.  If $A$ has the property that $a \precsim b$ whenever $s(a) < s(b)$
for every $s \in \mathrm{LDF}(A)$, then we say that $A$ has 
{\it strict comparison of positive elements}.

\section{A stability theorem for $W(\mathrm{C}(X))$}\label{stab}

\subsection{Strategy}

In this section we provide a positive answer to Question \ref{mainq}.  
The proof is long and proceeds in several steps, 
so an overview of our strategy is in order.  
We first generalise the concept of a trivial vector bundle by introducing
trivial positive elements in matrices over a commutative $C^*$-algebra,
and prove that our question need only be answered in the case that 
$X$ is a finite simplicial complex (Proposition \ref{reduce}). For a finite simplicial complex $X$, 
we prove that any positive element in $\mathrm{M}_n(\mathrm{C}(X))$
can be approximated in norm from below by a  particularly tractable
kind of positive element (Theorem \ref{wellsupapprox}).  Such an
approximant is then shown to be dominated (in the sense of Cuntz
comparison) by a trivial element whose rank over each point of $X$
exceeds the rank of the approximant by an amount no greater than
some fixed multiple of the dimension of $X$ (Theorem \ref{trivmajor}).
Finally, we prove that any positive element dominates 
a trivial element provided that the rank of the given element exceeds
the rank of the trivial element by $\mathrm{dim}(X)+1$ over every point
in $X$ (Theorem \ref{trivminor}).  Combining the last two theorems
yields a positive answer to Question \ref{mainq}.     

\begin{rems} {\rm
We are seeking the constant $K$ of Question \ref{mainq}.  The reader will
notice that we make little effort to find the smallest possible value
of $K$.  This is deliberate.  The applications
of Section \ref{ah} require only that some $K$ exists, and any effort to find
the optimal $K$ introduces more complexity into an already difficult proof.
Of course, it is interesting to ask what the optimal value of $K$ might be.
We suspect that it is the same as it is for projections, namely, about half
the covering dimension of $X$. }
\end{rems}

\subsection{Notation and background}

Let $X$ be a topological space.  By an upper (resp. lower) semicontinuous function 
$f:X \to \mathrm{M}_n(\mathbb{C})_+$ we will mean a function such that for every 
vector $\xi \in \mathbb{C}^n$, the real-valued function
\[
x \mapsto \langle f(x)\xi \ | \ \xi \rangle
\]
is upper (resp. lower) semicontinuous (cf. \cite{BE}).  The following result of Bratteli 
and Elliott, based on earlier work of D{\u{a}}d{\u{a}}rlat, Nagy, N{\'{e}}methi, and
Pasnicu (\cite[Proposition 3.2]{DNNP}), will be used extensively in the sequel.

\begin{thms}[Bratteli-Elliott, Theorem 3.1, \cite{BE}]
Let $X$ be a compact metrisable Hausdorff space of dimension $d$, and let $P$
and $Q$ be maps from $X$ into the projections of $\mathrm{M}_n$ such that
$P$ is lower semicontinuous and $Q$ is upper semicontinuous.  Suppose that
\[
P(x) \geq Q(x), \ \forall x \in X,
\]
and that, furthermore, there exists a natural number $k$ such that
\[
\mathrm{rank}(P(x)) > k + \frac{1}{2}(d+1), \ \forall x \in X,
\]
and 
\[
\mathrm{rank}(Q(x)) < k - \frac{1}{2}(d+1), \ \forall x \in X.
\]
Then, there is a continuous map $R$ from $X$ into the rank-k
projections of $\mathrm{M}_n$ such that
\[
P(x) \geq R(x) \geq Q(x), \ \forall x \in X.
\]
\end{thms}

For a natural number $l \in \mathbb{N}$, let $\theta_l$ denote the trivial vector
bundle of complex fibre dimension $l$.  In a metric space $X$ we will use the notation
$B_r(x)$ for the open ball of radius $r>0$ about $x \in X$.  The following 
well known theorem is a direct consequence of \cite[Theorems 1.2 and 1.5, Chapter 8]{H}.

\begin{thms}\label{projcomp}
Let $X$ be a compact metrisable Hausdorff space of dimension $d$, and 
let $\omega$ be a complex vector bundle over $X$.  Then, there is a complex vector
bundle $\overline{\omega}$ on $X$ of rank less than or equal to $d$ with the property
that $\omega \oplus \overline{\omega}$ is a trivial vector bundle.
\end{thms}

Let $\Gamma_n(X)$ denote the set of $n$-multisets (sets of $n$ elements, allowing
 multiplicity) whose entries are elements of a metric space $X$.  Equip $\Gamma_n(X)$ with the
following metric:  for $A,B \in \Gamma_n$,  let $\mathcal{P}$ denote the set of all
possible pairings of the elements of $A$ with the elements of $B$;  for $P \in \mathcal{P}$,
let $\Delta(P)$ denote the maximum distance between two paired elements in $P$;
set $\mathrm{dist}(A,B) = \mathrm{min}_{P \in \mathcal{P}} \Delta(P)$.

\subsection{Trivial positive elements}\label{trivpos}

\begin{dfs}\label{trivial}
Let $X$ be a compact metrisable Hausdorff space, and let $a \in \mathrm{M}_n(\mathrm{C}(X))$ 
be positive with (lower semicontinuous)
rank function $f:X \to \mathbb{Z}^+$ taking values in $\{n_1,\ldots,n_k\}$,
$n_1 < n_2 < \cdots < n_k$. \\ 
\begin{enumerate}
\item[(i)]
For each $1 \leq i \leq k$ define sets
\[
G_{i,a} := \{ x \in X | f(x) > n_i \},
\]
\[
F_{i,a} := \{ x \in X | f(x) = n_i \},
\] 
and
\[
H_{i,a} := \{ x \in X | f(x) \leq n_i \}.
\]
Note that $G_{i,a}$ is open. \\
\item[(ii)] Say that $a$ is trivial if there exist a natural number $n$, $n_k$ mutually orthogonal
projections $p_1,\ldots,p_{n_k} \in \mathrm{M}_{n}(\mathrm{C}(X))$, each corresponding
to a trivial complex line bundle,  and positive continuous functions 
$g_{i,a}:X \to \mathbb{R} \mathbf{1}_{\mathrm{M}_n(\mathrm{C}(X))}$ with
$\mathrm{supp}(g_{i,a}) = G_{i,a}$ such that $a$ is Cuntz equivalent to 
\[
\left( \bigoplus_{j=1}^{n_1} g_{1,a} p_j \right) \oplus \cdots
\oplus \left( \bigoplus_{j=n_{k-1}+1}^{n_k} g_{k,a} p_j \right).
\]
\\
\item[(iii)] Say that $a$ is well supported if, for each $1 \leq i \leq k$, there is a
projection $p_i \in \mathrm{M}_n(\mathrm{C}(\overline{F_{i,a}}))$ such that 
\[
\lim_{r \to \infty} a(x)^{1/r} = p_i(x), \ \forall x \in F_{i,a},
\]
and $p_i(x) \leq p_j(x)$ whenever $x \in \overline{F_{i,a}} \cap \overline{F_{j,a}}$ and
$i \leq j$.
\end{enumerate}
\end{dfs}

Clearly, if $a$ above is a projection, then it is trivial in the sense we have defined
if and only if it corresponds to a trivial vector bundle.  When there is no danger of confusion, we will 
omit the second subscript for the sets $G_{i,a}$, $F_{i,a}$, and $H_{i,a}$.


The following lemma was observed in \cite{Pe1}.  The proof is an easy exercise.

\begin{lms}\label{matrixcomp}
Let $a, b \in \mathrm{M}_n$ be positive.  Then, $a \precsim b$ if and only if
$\mathrm{rank}(a) \leq \mathrm{rank}(b)$.
\end{lms}




An important analogy between trivial positive elements and trivial vector bundles
is the fact that Cuntz comparison for trivial positive elements is encoded by rank 
functions.

\begin{props}\label{trivrankcomp}
Let $X$ be a compact Hausdorff space, and
let $a,b \in \mathrm{M}_n(\mathrm{C}(X))_+$ be trivial.  Then, 
$a \precsim b$ if and only if $\mathrm{rank}(a)(x) \leq \mathrm{rank}(b)(x)$,
$\forall x \in X$.
\end{props}

\begin{proof}
If $a \precsim b$, then $a(x) \precsim b(x)$ for every $x$ in $X$.  
It follows from Lemma \ref{matrixcomp} that 
\[
\mathrm{rank}(a)(x) \leq \mathrm{rank}(b)(x), \ \forall x \in X.
\]

Now suppose that the rank inequality above holds.
Let $n_1 < n_2 < \cdots < n_s$ (resp. $m_1 < m_2 < \cdots < m_t$) be the rank
values taken by $a$ (resp. $b$).
By assumption there exist projections $p_1,\ldots,p_{n_s}$ and $q_1,\ldots,q_{m_t}$ in
some $\mathrm{M}_n(\mathrm{C}(X))$, each corresponding to a trivial line bundle,
and positive continuous functions 
\[
g_{i,a}:X \to \mathbb{R} \mathbf{1}_{\mathrm{M}_n(\mathrm{C}(X))}, \ 1 \leq i \leq s,
\]
and
\[
g_{j,b}:X \to \mathbb{R} \mathbf{1}_{\mathrm{M}_n(\mathrm{C}(X))}, \ 1 \leq j \leq t,
\]
with $\mathrm{supp}(g_{i,a}) = G_{i,a}$ and $\mathrm{supp}(g_{j,b}) = G_{j,b}$
such that
\[
a \sim \left( \bigoplus_{j=1}^{n_1} g_{1,a} p_j \right) \oplus \cdots
\oplus \left( \bigoplus_{j=n_{s-1}+1}^{n_s} g_{s,a} p_j \right)
\] 
and
\[
b \sim \left( \bigoplus_{j=1}^{m_1} g_{1,b} q_j \right) \oplus \cdots
\oplus \left( \bigoplus_{j=m_{t-1}+1}^{m_t} g_{t,b} q_j \right).
\]
The projections $p_i$ and $q_i$ are Murray-von Neumann equivalent for 
every $1 \leq i \leq n_s$, as are their complements inside $\mathrm{M}_n(\mathrm{C}(X))$.
It follows that there is a unitary $U_1$ in $\mathrm{M}_n(\mathrm{C}(X))$ such that $u_1 q_1 u_1* = p_1$.
We can repeat this argument to find a unitary $u_2$ in $(1-p_1)\mathrm{M}_n(\mathrm{C}(X))(1-p_1)$
such that $u_2 q_2 u_2* = p_2$.  Then $U_2 = (u_2 \oplus p_1)u_1$ is a unitary in
$\mathrm{M}_n(\mathrm{C}(X))$ such that $U_2 q_i U_2^* = p_i$ for $i = 1,2$.
Iterating this process we arrive at a unitary $U$ in $\mathrm{M}_n(\mathrm{C}(X))$
such that
\[
U q_i U^* = p_i, \ 1 \leq i \leq n_s.
\]
Our rank inequality implies that 
\[
\mathrm{supp}(g_{i,a}) \subseteq \mathrm{supp}(g_{i,b}), \ 1 \leq i \leq n_s,
\]
whence $g_{i,a} p_i$ is in the hereditary subalgebra of $\mathrm{M}_n(\mathrm{C}(X))$
generated by $g_{i,b} p_i$.  It follows that 
\[
a \sim \left( \bigoplus_{j=1}^{n_1} g_{1,a} p_j \right) \oplus \cdots
\oplus \left( \bigoplus_{j=n_{s-1}+1}^{n_s} g_{s,a} p_j \right)
\] 
is in the hereditary subalgebra of $\mathrm{M}_n(\mathrm{C}(X))$ generated by
\[
U\left[ \left( \bigoplus_{j=1}^{m_1} g_{1,b} q_j \right) \oplus \cdots
\oplus \left( \bigoplus_{j=m_{t-1}+1}^{m_t} g_{t,b} q_j \right) \right]U^* \sim b,
\]
and that $a \precsim b$, as desired.
\end{proof}

\subsection{Reduction to finite simplicial complexes}
\begin{props}\label{reduce}
Suppose that there exists a constant $K > 0$ such that for any finite simplicial complex $X$ 
and positive elements $a,b \in \mathrm{M}_{\infty}(\mathrm{C}(X))$ for which
\[
\mathrm{rank}(a)(x) + K \mathrm{dim}(X) \leq \mathrm{rank}(b)(x), \ \forall x \in X,
\]
one has $a \precsim b$.  Then, this same statement holds (with the same value for $K$)
upon replacing finite simplicial complexes with arbitrary compact metric spaces.
\end{props}

\begin{proof}
Let $a$, $b$, and $K>0$ be as in the hypotheses of the theorem.
Let there be given a compact metric space $Z$, and positive
elements $a,b \in \mathrm{M}_{\infty}(\mathrm{C}(Z))$ such that
\[
\mathrm{rank}(a)(z) + K \mathrm{dim}(Z) \leq \mathrm{rank}(b)(z), \ \forall z \in Z.
\]
(This implies, in particular, that $Z$ has finite covering dimension, but
Theorem \ref{main} holds vacuously if $\mathrm{dim}(Z) = \infty$;  there is no loss of generality
here.)

A central theorem in the dimension theory of topological spaces asserts that if $Z$ is a 
compact metric space of finite covering dimension, then $Z$ is the limit of an inverse
system $(Y_i,\pi_{i,j})$, where each $Y_i$ is a finite simplicial complex of dimension
less than or equal to the dimension of $Z$ (cf. \cite[Theorem 1.13.2]{E}).  Thus, we have an inductive limit decomposition
for $\mathrm{M}_{\infty}(\mathrm{C}(Z))$:
\[
\mathrm{M}_{\infty}(\mathrm{C}(Z)) = \lim_{i \to \infty} \left( \mathrm{M}_{\infty}(\mathrm{C}(Y_i)),\phi_i \right),
\]
where $\phi_i:\mathrm{M}_{\infty}(\mathrm{C}(Y_i)) \to \mathrm{M}_{\infty}(\mathrm{C}(Y_{i+1}))$ is induced
by $\pi_{i,i+1}:Y_{i+1} \to Y_i$.

We claim that for each $\epsilon>0$ there exists a $\delta>0$ such that
\[
\mathrm{rank}(a-\epsilon)_+(z) + K \mathrm{dim}(Z) \leq \mathrm{rank}(b-\delta)_+(z), \ \forall z \in Z.
\]

Notice that if $\epsilon_1 \leq \epsilon_2$, then $\mathrm{rank}(a-\epsilon_1)_+ \geq \mathrm{rank}(a-\epsilon_2)_+$.
Let $\epsilon>0$ be given, and fix $z \in Z$.  Put $A = \mathrm{rank}(a)(z)$ and $B = 
\mathrm{rank}(b)(z)$ for convenience.  Let $\eta(z)$ denote the smallest non-zero eigenvalue
of $a(z)$.  For the purpose of proving our claim, we view $a$ and $b$ as being contained in
some $\mathrm{M}_n(\mathrm{C}(Z))$.   The map 
\[
\sigma:\mathrm{M}_n(\mathrm{C}(Z)) \times Z \to \Gamma_n
\] 
which assigns to an ordered pair $(d,x)$ the multiset whose elements are the eigenvalues of 
$d(x)$ is continuous in both variables, and so there is a neighbourhood
$V(z)$ of $z$ upon which $a$ has precisely $A$ eigenvalues greater than or equal to $\mathrm{min}\{
\eta(z)/2, \epsilon\}$.  It follows from the functional calculus that $\mathrm{rank}(a-\epsilon)_+$
is less than or equal to $A$ on $V(z)$.  
$\sigma(b,\bullet)$ is continuous, and so there is a neighbourhood $U(z)$ of $z$ upon which $b$ 
has at least $B$ non-zero eigenvalues.  In fact, more is true:  $U(z)$ may be chosen so that
there is a continuous choice of $B$ non-zero eigenvalues of $b$ on $U(z)$ which coincides with
the spectrum of $b(z)$ at $z$.  Let $\delta_z$ be half the smallest eigenvalue occuring in this
continous choice of eigenvalues.  Then, $\mathrm{rank}(b-\delta_z)_+ \geq B$ on $U(z)$, and
this remains true if $\delta_z$ is replaced by some smaller $\delta^{'}$.  Put $W(z) = V(z) 
\cap U(z)$.  Then,
\[
\mathrm{rank}(a-\epsilon)_+(x) \leq A +  K \mathrm{dim}(Z) \leq B \leq \mathrm{rank}(b-\delta_z)_+(x), \ \forall x \in W(z).
\]
$Z$ is compact, so we may find $z_1,\ldots,z_k \in Z$ such that 
\[
Z  = W(z_1) \cup \cdots \cup W(z_k).
\]
Put $\delta = \mathrm{min} \{\delta_{z_1},\ldots,\delta_{z_k}\}$, and let $x \in Z$ be given.
Then $x \in W(z_l)$ for some $1 \leq l \leq k$, and applying the preceding inequality we have
\[
\mathrm{rank}(a-\epsilon)_+(x) + K\mathrm{dim}(Z) \leq \mathrm{rank}(b-\delta_{z_l})_+(x) \leq
\mathrm{rank}(b-\delta)_+(x).
\]
This proves the claim.

Given $\eta>0$ and $a \in \mathrm{M}_{\infty}(\mathrm{C}(Z))_+$ one can always find an element
$\tilde{a}$ satisfying:
\begin{enumerate}
\item[(i)] $\tilde{a} \precsim a$;
\item[(ii)] $(a-\eta)_+ \precsim \tilde{a}$;
\item[(iii)] $\tilde{a}$ is in the image of $\phi_{i,\infty}$ for some $i \in \mathbb{N}$.
\end{enumerate}
(This follows from Proposition \ref{basics} and the inductive limit
decomposition of $\mathrm{M}_{\infty}(\mathrm{C}(Z))$.)  Let $\epsilon>0$ be given, and find 
$\delta>0$ so that
\[
\mathrm{rank}(a-\epsilon/2)_+(z) + K \mathrm{dim}(Z) \leq \mathrm{rank}(b-\delta)_+(z), \ \forall z \in Z.
\]
Find an approximant $\tilde{b}$ for $b$ satisfying (i)-(iii) above with $\eta=\delta$.  Similarly,
find an approximant $\tilde{a}$ for $(a-\epsilon/2)_+$ with $\eta=\epsilon/2$.  Then, $\forall z \in Z$,
\begin{eqnarray*}
\mathrm{rank}(\tilde{a})(z) + K\mathrm{dim}(Z) & \leq & \mathrm{rank}(a-\epsilon/2)_+(z) + K\mathrm{dim}(Z) \\
& \leq & \mathrm{rank}(b-\delta)_+(z) \\
& \leq & \mathrm{rank}(\tilde{b})(z).
\end{eqnarray*}

We may assume that both $\tilde{a}$ and $\tilde{b}$ are the images under $\phi_{i,\infty}$
of elements $\hat{a}$ and $\hat{b}$ in $\mathrm{M}_{\infty}(\mathrm{C}(Y_i))$, respectively. 
These pre-images satisfy
\[
\mathrm{rank}(\hat{a})(y) + K \mathrm{dim}(Z) \leq \mathrm{rank}(\hat{b})(y), \ \forall y \in \mathrm{Im}(\pi_{i,\infty}),
\]
where $\pi_{i,\infty}:Z \to Y_i$ is the continuous map which induces $\phi_{i,\infty}$. 
We cannot apply the hypothesis of the proposition unless the inequality above 
holds for all $y \in Y_i$, and so we modify the pre-image $\hat{b}$.  We may, as 
before, view $\hat{a}$ and $\hat{b}$ as lying in some $\mathrm{M}_n(\mathrm{C}(Y_i))$.
Choose a continuous function $f:Y_i \to [0,1]$ supported on the complement 
of $\mathrm{Im}(\pi_{i,\infty})$.  Put 
\[
\hat{\hat{b}} = \hat{b} \oplus f \cdot \mathbf{1}_{\lceil K\mathrm{dim}(Z) \rceil}.
\]
Since $\hat{\hat{b}}$ and $\hat{b}$ agree on $\mathrm{Im}(\pi_{i,\infty})$, we have
that $\phi_{i,\infty}(\hat{\hat{b}}) = \tilde{b}$.  But clearly
\[
\mathrm{rank}(\hat{a})(y) + K \mathrm{dim}(Z) \leq \mathrm{rank}(\hat{b})(y), \ \forall y \in Y_i,
\]
whence $\hat{a} \precsim \hat{\hat{b}}$ by our hypothesis.  Cuntz subequivalence is preserved
under $*$-homomorphisms, whence $\tilde{a} \precsim \tilde{b}$.  Now
\[
(a-\epsilon)_+ \precsim \tilde{a} \precsim \tilde{b} \precsim b;
\]
$\epsilon$ was arbitrary, and the proposition follows. 
\end{proof}

\subsection{Well supported approximants}
\begin{lms}\label{cellapprox}
Let $X$ be a finite simplicial complex, $V$ be an open subset of $X$, and $U$ be
a closed subset of $V$.  Then, there is a refinement of 
the simplicial structure on $X$ and a subcomplex $Y$ of this refined structure
satisfying:
\begin{enumerate}
\item[(i)] $Y \supseteq V^c$;
\item[(ii)] $U \cap Y = \emptyset$;
\item[(iii)] $Y$ is the closure of its interior.
\end{enumerate}  
Moreover, $\overline{Y^c}$ and $\partial Y = \partial Y^c$ are subcomplexes of this 
refined structure.
\end{lms}

\begin{proof}
We first define a precursor $\widetilde{Y}$ to $Y$, whose definition we later
refine to obtain $Y$ proper.
By assumption, $U \cap V^c = \emptyset$.  Since $U$ and
$V^c$ are compact, there is a $\delta > 0$ such that 
\[
\mathrm{dist}(U,x) > \delta, \ \forall x \in V^c.
\]
Refine the simplicial structure on $X$ through repeated barycentric subdivision
until the largest diameter of any simplex is less than $\delta/2$.   Let $\widetilde{Y}$ be 
the subcomplex consisting of all simplices whose intersection with $V^c$ is
nonempty.  The distance from any point in $\widetilde{Y}$ to $V^c$ is at most $\delta/2$, and
so $\widetilde{Y} \cap U = \emptyset$.  Every $x \in V^c$ is contained in some simplex
of $X$, so $\widetilde{Y} \supseteq V^c$. 

Now choose $x \in \widetilde{Y}^c$, and let $\Theta_x$ be the smallest simplex of the refined simplicial
structure on $X$ containing
$x$.  Suppose that $\Theta_x$ is not contained in $\overline{\widetilde{Y}^c}$.  Then, $\Theta_x$
contains a point $y \in \widetilde{Y}^{\circ}$, and there is an open set $U \subseteq \widetilde{Y}^{\circ}$
such that $y \in U$.  Thus, there is a point $y^{'} \in U$ which is in the (relative) interior
of $\Theta_x$, and the smallest simplex of $X$ containing $y^{'}$ is $\Theta_x$.  This implies that $\Theta_x$ is contained
in $\widetilde{Y}$ by construction, contradicting $x \in \Theta_x$.  We conclude that
$\Theta_x \subseteq \overline{\widetilde{Y}^c}$.  Now 
\[
\widetilde{Y}^c \subseteq \bigcup_{x \in \widetilde{Y}^c} \Theta_x \subseteq \overline{\widetilde{Y}^c},
\]
and so the second containment above is in fact equality.  We conclude
that $\overline{\widetilde{Y}^c}$ is a subcomplex of $X$.  $\widetilde{Y}^{c}$
is open, whence $\overline{\widetilde{Y}^c}$ is the closure of its interior.

Notice that $\overline{\widetilde{Y}^c}$ is defined in the same manner as $Y$:
for each point in a given open set, one finds the smallest simplex containing the 
said point, and then takes the union of these simplices over all points in the
open set.  We may therefore define 
\[
Y := \overline{\left( \overline{\widetilde{Y}^c} \right)^c},
\]
and apply the arguments above to conclude that $Y$ is a subcomplex.  $Y$ both
contains $V$ and is contained in $\widetilde{Y}$, and so satisfies conclusions (i)
and (ii) of the Lemma;  $Y$ satisfies conclusion (iii) by construction.  Repeating
the arguments above one last time, we find that $\overline{Y^c}$ and hence
$\partial Y^c = \partial Y$ are subcomplexes.
\end{proof}

\begin{thms}\label{wellsupapprox}
Let $X$ be a finite simplicial complex, $a \in \mathrm{M}_n(\mathrm{C}(X))_+$, and $\epsilon > 0$
be given.  Then, there exists a well supported element
$f \in \mathrm{M}_n(\mathrm{C}(X))_+$ such that $f \leq a$
and $||f-a|| < \epsilon$.  Moreover, $f$ takes the same 
rank values as $a$, and the sets $\overline{F_{i}}$ corresponding to $f$ (see Definition \ref{trivial}) 
may be assumed, upon refining the simplicial structure 
of $X$, to be subcomplexes of $X$.
\end{thms}

\begin{proof}
Let $a \in \mathrm{M}_n(\mathrm{C}(X))_+$, and let $H_i$ be the 
set $H_{i,a}$ of Definition \ref{trivial}.  Let $\epsilon > 0$ be given.  

\vspace{2mm}
\noindent
{\it Step 1.} We will show that there exist open sets
$V_i \supseteq H_i$, $1 \leq i \leq k$ and an upper semicontinuous
function $g:X \to \mathrm{M}_n(\mathbb{C})_+$ satisfying:
\begin{enumerate}
\item[(i)] $g(x) = a(x)$, $\forall x \in H_i \backslash (\cup_{j=1}^{i-1} V_j)$;
\item[(ii)] $g(x)$ is continuous on  $V_i \backslash (\cup_{j=1}^{i-1} V_j)$, $1 \leq i \leq k$;
\item[(iii)] $||g(x)-a(x)|| < \epsilon$, $\forall x \in X$;
\item[(iv)] $\mathrm{rank}(g)(x) = n_i$, $\forall x \in V_i \backslash (\cup_{j=1}^{i-1} V_j)$;
\item[(v)] with $D_i := \overline{V_i} \backslash (\cup_{j=1}^{i-1} V_j)$,
there is a continuous projection-valued function
$p_i \in \mathrm{M}_n(\mathrm{C}(D_i))$ such that
\[
\lim_{r \to \infty} g(x)^{1/r} = p_i(x), \ \forall x \in V_i \backslash (\cup_{j=1}^{i-1} V_j),
\]
and $p_i(x) \leq p_j(x)$ whenever $x \in D_i \cap D_j$ and
$i \leq j$.
\end{enumerate}
Moreover, upon refining the simplicial structure of $X$, we may assume that both
$V_i$ and $\overline{V_i^c}$ are subcomplexes of $X$.

Let $n_1 < \cdots < n_k$ be the rank values taken by $a$.  The function which
assigns to each point $x \in \partial H_1$ the minimum non-zero eigenvalue of $a(x)$ 
is continuous on a compact set, and so achieves a minimum, say $\eta_1 > 0$.  For each $x \in \partial H_1$, find
$\delta_x > 0$ such that for every $y \in B_{\delta_x}(x)$ one has that the
eigenvalues of $a(y)$ are all either greater than 
\[
L_1^u := \mathrm{max}\{(2/3)\eta_1, \eta_1-\epsilon \},
\]
or less than
\[
L_1^l := \mathrm{min}\{(1/3)\eta_1, \epsilon\}.  
\]
Let $U_x$ be the connected component of $B_{\delta_x}(x)$ containing $x$.
Put 
\[
W_1 = H_1 \cup \left( \bigcup_{x \in \partial H_1} U_x \right).
\]
$W_1$ is open.  
Use Lemma \ref{cellapprox} to find
a subcomplex $Y_1$ of $X$ such that $Y_1 \supseteq \widetilde{W_1}^c$ and
$Y_1 \cap H_1 = \emptyset$.  Put $V_1 = Y_1^c$, and note that $V_1$ so
defined is a subcomplex of some (possibly refined) simplicial structure on $X$. 
Define $p_1(x)$ to be the support projection of $a(x)$ for $x \in H_1$, and
the support projection of the eigenvectors of $a(x)$ corresponding to eigenvalues 
which are greater than or equal to $L_1^u$ for $x \in \overline{V_1} \backslash H_1$.
For each $x \in V_1$, put $g(x) = p_1(x) a(x)$.  

\vspace{2mm}
\noindent
{\it Claim:}
\begin{enumerate}
\item[(i)] $g(x) = a(x)$, $\forall x \in H_1$;
\item[(ii)] $g(x)$ is continuous on $V_1$;
\item[(iii)] $||g(x)-a(x)|| < \epsilon$, $\forall x \in V_1$;
\item[(iv)] $\mathrm{rank}(g)(x) = n_1$, $\forall x \in V_1$;
\item[(v)] $p_1$ is a continuous
projection valued function on $\overline{V_1}$ such that
\[
\lim_{r \to \infty} g(x)^{1/r} = p_i(x), \ \forall x \in V_1.
\]
\end{enumerate}

\vspace{2mm}
\begin{proof} (i) is clear from the definition of $g$ on $H_1$.  

For (ii), let
$x_n \to x$ in $V_1$.  If $x$ is an interior point of $H_1$, then $g(x_n) \to g(x)$ since
$a$ is continuous and $g = a$ on $H_1$.  Otherwise, $x \in V_1$ is an interior point of some
$U_y$, and we may find $N \in \mathbb{N}$ such that $x_n \in U_y$, $\forall n \geq N$.
It will thus suffice to prove that $g(x)$ is continuous on $U_y$, $y \in \partial H_1$.
Let $\sigma(a(x))$ denote the spectrum of $a$ at the point $x \in X$.  
  The map $s:X \to \Gamma_n$ given by
$s(x) = \sigma(a(x))$ is continuous.  Thus, $s(x_n) \to s(x)$.  In particular, 
the submultiset of $s(x_n)$ corresponding to elements larger than $L_1^u$ 
converges to the similar submultiset of $s(x)$;  the eigenvectors
of $a(x_n)$ corresponding to $p_1(x_n)$ converge to the eigenvectors if $a(x)$ corresponding
to $p(x)$, and so $g(x_n) \to g(x)$, as required.  

(iii) follows from the functional calculus and the definition of $g$ ---
$a(x)-g(x)$ is zero if $x \in H_1$, and is equal, in the functional
calculus, to $h(t)=t$ on the part of the spectrum of $a$ which is less than
$L_1^l \leq \epsilon$ and zero otherwise.

(iv) is trivial for $x \in H_1$, so suppose that $x \in U_y$ for 
some $y \in \partial H_1$.  The rank of $g(x)$ is equal to the number of eigenvalues of $a(x)$
which are greater than $L_1^u$, and this number is $n_1$
for $a(y)$.  Let $M \subseteq \Gamma_n$ consist of those multisets with the property that
each element of the multiset is either greater than $L_1^u$
or less than $L_1^l$.  Let $A \subseteq M$ consist of those
multisets for which the number of elements greater than $L_1^u$
is exactly $n_1$, and let $B$ be the complement of $A$ relative to $M$.  Then, $A$ and $B$ are
separated.  Since $s(z) \subseteq M$, $\forall z \in U_y$, $s(y) \in A$, and $U_y$ is connected,
we conclude that $s(z) \in A$, $\forall z \in U_y$.  In particular, the rank of $g(x)$ is $n_1$,
as required.  

Observe that (iv) implies that the function which assigns to a point $x \in V_1$
the submultiset of $s(x)$ consisting of those eigenvalues which are greater than or
equal to $L_1^u$ is continuous on $V_1$.  $p_1(x)$ is the spectral projection on
this submultiset for each $x \in V_1$, and is thus continuous merely by the existence
of the continuous functional calculus.  This proves (v), and hence the claim.
\end{proof}

The claim above (now proved) is the base case of an inductive argument.
We now describe the construction of $V_2$, $p_2$, and the definition of $g$ 
on $V_2 \backslash V_1$.  The construction of the subsequent $V_i$s, $p_i$s, and the
definition of $g$ on $V_i \backslash (\cup_{j=1}^{i-1} V_j)$ will be
similar --- all of the essential difficulties are encountered already for $i=2$.

Let $\partial (H_2 \backslash V_1)$
be the boundary of $H_2 \backslash V_1$ inside $V_1^c$ (we may assume that this boundary
is nonempty by shrinking $V_1$, if necessary). Find, as before,
the minimum value $\eta_2 >0$ occuring as a non-zero eigenvalue of $a|_{\partial(H_2 \backslash V_1)}$.
Define
\[
L_2^u = \mathrm{max}\{(2/3)\eta_2,\eta_2-\epsilon\},
\]
and
\[
L_2^l = \mathrm{min}\{(1/3)\eta_2,(1/3)\eta_1,\epsilon\}.
\]
(Note the dependence of $L_2^l$ on $\eta_1$.)
Find, for each $x \in \partial (H_2 \backslash V_1)$, a connected open (rel $V_1^c$) set $U_x \subseteq V_1^c$ containing
$x$ and with the property that for every $y \in U_x$, the eigenvalues of $a(y)$ are all either
greater than $L_2^u$ or less than $L_2^l$.  Let $\widetilde{U_x}$ be an open set in $X$ such that 
$U_x = \widetilde{U_x} \cap V_1^c$.  Put 
\[
W_2 = H_2 \cup \left(\bigcup_{x \in \partial(H_2 \backslash V_1)} \widetilde{U_x} \right) \subseteq X. 
\]
Refine the simplicial structure on $X$ so that we may find (using Lemma \ref{cellapprox}) a subcomplex $Y_2$ of $V_1^c$ which
is disjoint from $W_2^c$ and whose interior (rel $V_1^c$) contains $H_2 \backslash V_1$.  Put $V_2 = Y_2^{\circ} \cup V_1$ (interior
taken (rel $V_{1}^c$)).  Thus, $V_2^c$ 
and $\overline{V_2}$ are subcomplexes of $X$.
Define $p_2(x)$ to be the support projection of $a(x)$ for $x \in H_2 \backslash V_1$.  For
any $y \in (U_x \backslash H_2) \cap V_2$, some $x \in \partial(H_2 \backslash V_1)$, let $p_2(x)$ be the support
projection of those eigenvectors of $a(y)$ corresponding to eigenvalues greater than
$L_2^u$.  Put $g(x) = p_2(x) a(x)$, $\forall x \in V_2 \backslash V_1$.  

\vspace{2mm}
\noindent
{\it Claim:}  $p_2$ and $g(x)$ so defined have the desired properties, namely,
\begin{enumerate}
\item[(i)] $g(x) = a(x)$, $\forall x \in H_2 \backslash V_1$;
\item[(ii)] $g(x)$ is continuous on $V_2 \backslash V_1$;
\item[(iii)] $||g(x)-a(x)|| < \epsilon$, $\forall x \in V_2$;
\item[(iv)] $\mathrm{rank}(g)(x) = n_2$, $\forall x \in V_2 \backslash V_1$;
\item[(v)] $p_2$ extends to a continuous
projection valued function on $D_2 := \overline{V_2} \backslash V_1$ such that
\[
\lim_{r \to \infty} g(x)^{1/r} = p_i(x), \ \forall x \in V_1,
\]
and $p_1(x) \leq p_2(x)$ whenever $x \in \overline{V_1} \cap (\overline{V_2} \backslash V_1)$.
\end{enumerate}

\vspace{2mm}
\begin{proof}
The proofs of all but the last part of (v) are identical to the proofs of
the corresponding statements for $p_1$ above.  We must show that  $p_1(x) \leq p_2(x)$ whenever $x \in 
\overline{V_1} \cap (\overline{V_2} \backslash V_1)$.  This follows from the fact that the eigenvalues of 
$(1-p_2(x))a(x)(1-p_2(x))$ are all less than or equal to $L_2^l \leq L_1^l$, and so
correspond to eigenvectors in the complement of the range of $p_1(x)$.
\end{proof}

The remaining $p_i$s, $i < k$, may be constructed inductively in a manner similar to the construction
of $p_2$:  let $\eta_i$ be the minimum eigenvalue taken by $a$ on $\partial (H_i \backslash
 (\cup_{j=1}^{i-1} V_j))$ --- boundary relative to $X \backslash (\cup_{j=1}^{i-1} V_j)$ both
here and below, and assumed to be nonempty by shrinking $V_1,\ldots,V_{i-1}$ if necessary --- and put
\[
L_i^u = \mathrm{max}\{(2/3)\eta_i,\eta_i-\epsilon\},
\]
and 
\[
L_i^l = \mathrm{min}\{(1/3)\eta_i, L_{i-1}^l\};
\]
find connected open (rel $X \backslash  (\cup_{j=1}^{i-1} V_j)$) sets $U_x$ for each
$x \in \partial (H_i \backslash (\cup_{j=1}^{i-1} V_j)$ which contain
$x$, and have the property that every eigenvalue of $a(y)$, $y \in U_x$, is either
greater than $L_i^u$ or less than $L_i^l$;  find open sets $\widetilde{U_x} \subseteq
X$ such that $U_x = \widetilde{U_x} \cap X \backslash (\cup_{j=1}^{i-1} V_j)$, and
define 
\[
W_i = H_i \cup \left(\bigcup_{x \in \partial(H_i \backslash (\cup_{j=1}^{i-1} V_j) )} \widetilde{U_x} \right) \subseteq X; 
\]
refine the simplicial structure on $X$ and find (using Lemma \ref{cellapprox}) a subcomplex $Y_i$ of $V_{i-1}^c$ which
is disjoint from $W_i^c$ and whose interior (rel $V_{i-1}^c$) contains $H_i \backslash V_{i-1}$; put $V_i = Y_i^{\circ} \cup V_{i-1}$, (interior
taken (rel $V_{i-1}^c$)) so that $V_i^c$ and $\overline{V_i}$ are subcomplexes of $X$;
define $p_i(x)$ to be the support projection
of the eigenvectors of $a(x)$ with non-zero eigenvalues for $x \in H_i \backslash (\cup_{j=1}^{i-1} V_j)$,
and the support projection of those eigenvectors of $a(x)$ having eigenvalues greater than or
equal to $L_i^u$ for $x \in \overline{V_i} \backslash (\cup_{j=1}^{i-1} V_j)$;  put
$g = p_i(x) a(x)$, $\forall x \in V_i \backslash  (\cup_{j=1}^{i-1} V_j)$.

$p_i$ and $g$ so defined have the desired properties.
As before, all but the last part of statement (v) follow from the proofs of the
corresponding facts for $p_1$ and $p_2$.  So suppose that $x \in D_i \cap D_j$,
where $j \leq i$.  The range of $1-p_i(x)$ correponds to the span of eigenvectors of
$a(x)$ with eigenvalues less than or equal to $L_i^l \leq L_j^l$.  This range
is contained in the range of $1-p_j(x)$, since the latter corresponds to the span of 
eigenvectors of $a(x)$ with eigenvalues less than $L_j^l$.  It follows that $p_j(x)
\leq p_i(x)$, as required.   

For $i=k$, the situation is straightforward.  Simply put $g(x) = a(x)$,
$\forall x \in V_{k-1}^c$.

\vspace{2mm}
\noindent
{\it Step 2.}
We have $V_1 \subseteq \cdots \subseteq V_k = X$.  Since the $H_i$ are closed, we may find open
sets $U_1 \subseteq \cdots \subseteq U_k = X$ such that $H_i \subseteq U_i \subseteq \overline{U_i} 
\subseteq V_i$.  Moreover, we may assume, using Lemma \ref{cellapprox}, that $U_i^c$ and $\overline{U_i}$ are subcomplexes
(possibly empty) of $X$ for $1 \leq i \leq k$.
Let us define a positive upper semicontinuous function $\tilde{f}:X \to
\mathrm{M}_n(\mathrm{C}(X))$ as follows:  $\tilde{f}(x) = p_i(x) a(x)$, $\forall x \in
U_i \backslash (\cup_{j=1}^{i-1} U_j)$.  Then, $\tilde{f}(x) \geq g$, $\forall
x \in X$.  Moreover, $\tilde{f}$ is well supported by the $p_i$s.  To 
see this, one need only check the coherence condition.  Put $E_i = \overline{U_i} \backslash
(\cup_{j=1}^{i-1} U_j)$, and let $x \in E_i \cap E_j$, $j \leq i$.  The range of $1-p_i(x)$ correponds to the span of eigenvectors of
$a(x)$ with eigenvalues less than or equal to $L_i^l \leq L_j^l$ (since $U_i \subseteq V_i$).  This range
is contained in the range of $1-p_j(x)$, since the latter corresponds to the span of 
eigenvectors of $a(x)$ with eigenvalues less than $L_j^l$.  It follows that $p_j(x)
\leq p_i(x)$, as required.     

We now describe a smoothing process which will transform $\tilde{f}$ into
the function $f$ required by the theorem.  For the first step in our process 
we work inside $X \backslash U_{k-2}$.  Find a continuous function $s: X \backslash U_{k-2}
\to [0,1]$ which is zero on $\overline{U_{k-1}} \backslash U_{k-2}$, one on
$X \backslash V_{k-1}$, and non-zero off $\overline{U_{k-1}} \backslash U_{k-2}$.  Define $f^{k-1}$ on 
$X \backslash U_{k-2}$ as follows:
\[
f^{k-2}(x) = \left\{ \begin{array}{ll} \tilde{f}(x), & x \in \overline{U_{k-1}} \backslash U_{k-2}
\cup (X \backslash V_{k-1}) \\
g(x) \oplus s(x)(\tilde{f} \ominus g)(x), & x \in V_{k-1} \backslash
\overline{U_{k-1}} 
\end{array} \right. .
\]
$f^{k-2}$ is continuous on $X \backslash U_{k-2}$ and dominates $g$ on $X \backslash U_{k-2}$.
It is also subordinate to $a$.

For the generic step in our process, we work inside $X \backslash U_{k-i}$, $i \geq 3$.  Assume that we have
found a continuous function $f^{k-i+1}: X \backslash U_{k-i+1} \to \mathrm{M}_n(\mathbb{C})_+$
such that 
\[
g(x) \leq f^{k-i+1}(x), \ \forall x \in (X \backslash U_{k-i+1})
\cap V_{k-i+1}.
\]
(This is the key property required to go from one stage in the smoothing process to the 
next.)  Find a continuous function $s: X \backslash U_{k-i} \to [0,1]$ which is zero on
$\overline{U_{k-i+1}} \backslash U_{k-i}$, one on $X \backslash V_{k-i+1}$, and non-zero off
$\overline{U_{k-i+1}} \backslash U_{k-i}$.  Define
$f^{k-i}$ on $X \backslash U_{k-i}$ as follows:
\[
f^{k-i}(x) = \left\{ \begin{array}{ll} \tilde{f}(x), & x \in \overline{U_{k-i+1}} \backslash U_{k-i}\\
f^{k-i+1}(x), & x \in (X \backslash V_{k-i+1}) \\
g(x) \oplus s(x)(\tilde{f}^{k-i+1} \ominus g)(x), & x \in V_{k-i+1} \backslash
\overline{U_{k-i+1}} 
\end{array} \right. .
\]
Then, as before, $f^{k-i}$ is continuous on $X \backslash U_{k-i}$, dominates
$g$, and is subordinate to $a$.  This smoothing process terminates when $i = k-1$,
and the resulting continuous function $f:X \to \mathrm{M}_n(\mathbb{C})_+$ has
the desired properties.  (Note that the sets $H_i$ of Definition \ref{trivial} corresponding
to $f$ are precisely the $\overline{U_i}$.  The sets $\overline{F_i}$ of Definition \ref{trivial}
corresponding to $f$ are thus the $\overline{U_i} \cap U_{i-1}^c$, and so are
subcomplexes of the refined simplicial structure on $X$.)  
\end{proof}

\subsection{Trivial majorants}
\begin{lms}\label{extendone}
Let $X$ be a finite simplicial complex of dimension $d$, and let $Y$ be a subcomplex.  Let $p
\in \mathrm{M}_n(\mathrm{C}(Y))$ be a projection of (complex) fibre dimension $l \in \mathbb{N}$
corresponding to a trivial vector bundle, and suppose that $l \geq \lceil d/2 \rceil +1$.
Then, there is a projection $q \in A := \mathrm{M}_n(\mathrm{C}(X))$ such that
\begin{enumerate}
\item[(i)] $q$ corresponds to a trivial vector bundle of fibre dimension $l$;
\item[(ii)] $q(y) = p(y)$, $\forall y \in Y$.
\end{enumerate}   
\end{lms}

\begin{proof}
The projection $p$ may be viewed as a vector bundle $(E,r,Y)$. $E$ is trivial, and so admits $l$ mutually orthogonal
and everywhere non-zero cross-sections $s_i:Y \to E$, $1 \leq i \leq l$.  We will prove that each $s_i$ can be extended to
a continuous map $v_i$ defined on all of $X$ which takes values in $\mathbb{C}^n \backslash \{0\}$, and that these extensions
can be chosen to be mutually orthogonal.  The projection whose range at a point $x \in X$ is $\mathrm{span}\{v_1(x),\ldots,
v_l(x)\}$ is then the projection $q$ that we seek.

We proceed by induction.  First consider $s_1$, which may be viewed as a continuous map from
$Y$ to $\mathbb{C}^n \backslash \{0\} \cong \mathbb{R}^{2n} \backslash \{0\}$.  Theorem 2.2 
in Chapter 1 of \cite{H} states that if $(A,B)$ is a relative CW-complex and $R$ is a space 
which is connected in each dimension for which $A$ has cells, then every continuous map 
$f:B \to R$ extends to a continuous map $g:A \to R$.  $(X,Y)$ is a relative CW-complex, 
$\mathbb{C}^n$ is $(2n-1)$-connected, and $2n-1 \geq d$, so $s_1$ extends to a continuous
map $v_1:X \to \mathbb{C}^n \backslash \{0\}$, as desired.

Suppose now that we have found mutually orthogonal and continous extensions $v_i$ of $s_i$
for each $i < k \leq l$.  We wish to extend $s_k$ to a continous map $v_k$ on $X$ taking
values in $E$ and pointwise orthogonal to $v_1,\ldots,v_{k-1}$.  Let $Q_{k-1} \in
\mathrm{M}_n(\mathrm{C}(X))$ be the projection whose range at a point $x \in X$ is
$\mathrm{span}\{v_1(x),\ldots,v_{k-1}(x)\}$.  To find our extension $v_k$, it will
suffice to extend $s_k$ to an everywhere non-zero cross-section inside
\[
A_{k-1} := (\mathbf{1}_A - Q_{k-1}) (\mathrm{M}_n(\mathrm{C}(X))) (\mathbf{1}_A - Q_{k-1}) \cong \mathrm{M}_{n-k}(\mathrm{C}(X))      
\]
(the last isomorphism follows from the fact that $\mathbf{1}_A - Q_{k-1}$ is trivial).
This problem is identical to the problem of extending $s_1$, save that we are now
extending a map from $Y$ into $\mathbb{C}^{n-k} \backslash \{0\}$ rather than
into $\mathbb{C}^n \backslash \{0\}$.  Since $n-l \geq \lceil d/2 \rceil+1$, 
we have that $2(n-k)-1 \geq d$.  Since $\mathbb{C}^{n-k} \backslash \{0\}$
is $2(n-k)-1$-connected, we may use Theorem 2.2, Chapter 1, \cite{H}, to find 
the desired extension $v_k$.  

\end{proof}

\begin{lms}\label{trivcomp}
Let $X$ be a finite simplicial complex of dimension $d$, and let $U_1, \ldots, U_k$ 
be subcomplexes.  Suppose that for each $1 \leq i \leq k$ there is a constant rank
projection-valued 
map $p_i \in \mathrm{M}_{\infty}(\mathrm{C}(U_i))$ satisfying:
\begin{enumerate}
\item[(i)] $p_i(x) \leq p_j(x)$ whenever $x \in U_i \cap U_j$ and $i \leq j$; 
\item[(ii)] $\mathrm{rank}(p_i) < \mathrm{rank}(p_j)$ whenever $i < j$.
\end{enumerate}  

\noindent
Then, there exist a partition
of $\{1,2,\ldots,k\}$ into nonempty subsets $J_1,\ldots,J_s$, subcomplexes
$V_l := \cup_{i \in J_l} U_i$, and projections $R_l \in \mathrm{M}_{\infty}(\mathrm{C}(V_l))$,
$1 \leq l \leq s$, such that:
\begin{enumerate}
\item[(i)] each $R_l$ corresponds to a trivial constant rank vector bundle on $V_l$;
\item[(ii)] $3d+3 \leq \mathrm{rank}(R_l) - \mathrm{rank}(p_i) < 4d + 3$, $\forall i \in J_l$;
\item[(iii)] $\mathrm{rank}(R_l)+ d \leq \mathrm{rank}(R_{l+1})$, $1 \leq l < s$;
\item[(iv)] $p_i(x) \leq R_l(x)$ for each $x \in U_i \cap V_l$;
\item[(v)] if $x \in V_l \cap V_t$ and $l \leq t$, then $R_l(x) \leq R_t(x)$.
\end{enumerate}
\end{lms}

\begin{proof}
Put $n_i = \mathrm{rank}(p_i)$.  For each $n \in \mathbb{N}$, put 
\[
D_n = \{m \in \mathbb{N}| (n-1)d \leq m < nd \}.
\]
Let $\widetilde{J_1},\ldots,\widetilde{J_s}$ be the list of $D_n$s which have nonempty intersection with 
$\{n_1,\ldots,n_k\}$, ordered so that some (and hence every) element of $\widetilde{J_l}$ is less than every
element of $\widetilde{J_{l+1}}$.  For $1 \leq j \leq s$, let $J_j$ be the set of
indices of the $n_i$s appearing in $\widetilde{J_j}$.  

For $A \subseteq X$ put
\[
P_A(x) := \bigvee_{1 \leq i \leq k} p_i(x), \ \forall x \in A.
\]
The condition that  $p_i(x) \leq p_j(x)$
whenever $x \in U_i \cap U_j$ and $i \leq j$ implies that $P_A(x)$ is an 
upper semicontinuous projection-valued function for every $A \subseteq X$.

An application of \cite[Theorem 3.1]{BE} allows us to find a constant rank 
projection-valued map $Q_s \in \mathrm{M}_{\infty}(\mathrm{C}(V_s))$ such that
\[
\mathrm{rank}(Q_s) - \mathrm{rank}(p_k) \leq d
\]
and
\[
p_i(x) \leq Q_s(x), \ \forall i \in J_s, \ \forall x \in U_i.
\]
It follows that
$\mathrm{rank}(Q_s) - \mathrm{rank}(r_i) < 2d$ for every $i \in J_s$.  
For every $N \geq d$, there is a projection $\overline{Q_{s,N}} 
\in \mathrm{M}_{\infty}(\mathrm{C}(V_s))$ of rank $n$ such that $Q_s \oplus \overline{Q_{s,N}}$
corresponds to a trivial vector bundle (cf. Theorem \ref{projcomp}).  Put 
\[
N_s = 3d+3 + \mathrm{rank}(r_k) - \mathrm{rank}(Q_s)
\]
and
\[
R_s = Q_s \oplus \overline{Q_{s,N_s}}.
\]
$R_s$ so chosen satisfies conditions (i), (ii), and (iv) in the conclusion of the lemma;
conditions (iii) and (v) are not yet relevant.

Now suppose that we have found, for each $m < l \leq s$, constant rank projections
$R_l \in \mathrm{M}_{\infty}(\mathrm{C}(V_l))$ satisfying conditions (i)-(iv) of
the conclusion of the lemma, and satisfying condition (v) whenever $t,l > m$. 
We will construct $R_m$ on $V_m$ so that $R_1,\ldots,R_m$ satisfy (i)-(iv), and
satisfy (v) when $t,l \geq m$.  Proceeding inductively then yields the lemma.  

Define a projection-valued map 
\[
\widetilde{R}_{m}:\bigcup_{m \leq l \leq s} V_l \to \mathrm{M}_{\infty}(\mathbb{C})
\]
by setting $\widetilde{R}_{m}(x) = R_l(x)$ if $l >m$ is the smallest index such that
$x \in V_l$, and setting $\widetilde{R}_{m}(x)$ equal to the unit of the (arbitrarily large)
matrix algebra which constitutes the target space of all of our projection-valued maps
otherwise.  One easily checks that $\widetilde{R}_{m}$ is lower semicontinuous, and that
\[
\mathrm{rank}(\widetilde{R}_m-P_{V_m})(x) \geq 3d+3, \ \forall x \in V_m.
\]
An application of 
\cite[Theorem 3.1]{BE}  then yields a constant rank projection-valued
map $Q_m \in \mathrm{M}_{\infty}(\mathrm{C}(V_m))$ such that 
\[
P_{V_m}(x) \leq Q_m(x) \leq \widetilde{R}_{m}(x), \ \forall x \in V_m,
\]
and
\[
\mathrm{rank}(\widetilde{R}_m-Q_m)(x) \geq 2d+2, \ \forall x \in V_m.
\]
Applying \cite[Theorem 3.1]{BE} to $\widetilde{R}_m$ and $Q_m$ yields
a constant rank projection-valued map $Q_m^{'} \in \mathrm{M}_{\infty}(\mathrm{C}(V_m))$
such that
\[
Q_m(x) \leq Q_m^{'}(x) \leq \widetilde{R}_{m}(x), \ \forall x \in V_m,
\]
and
\[
\mathrm{rank}(Q_m^{'}-Q_m)(x) \geq d+1, \ \forall x \in V_m.
\]
By \cite[Chapter 8, Theorem 1.2]{H}, there is a  
subprojection $\overline{Q_m}$ of $Q_m^{'}-Q_m$ with constant rank $\leq d$ 
such that $T_m := Q_m \oplus \overline{Q_m}$
corresponds to a trivial vector bundle over $V_m$.  By definition we have
\[
P_{V_m}(x) \leq T_m(x) \leq \widetilde{R}_m(x), \ \forall x \in V_m.
\]

From the definition of $\widetilde{R}_m(x)$ we have
\[
P_{V_m}(x) \leq T_m(x) \leq \bigvee_{m < l \leq s} R_l(x), \ \forall x \in V_m \cap (V_{m+l} \cup \cdots \cup V_s).
\]     
$R_{m+1}-T_m$ is a trivial projection on $V_m \cap V_{m+1}$, and so can be extended
to a trivial projection of the same rank on $V_m \cap (V_{m+1} \cup V_{m+2})$ by
Lemma \ref{extendone} (this requires the fact that the $R_j$, $m < j \leq s$, satisfy
condition (iii) in the conclusion of the lemma).  We can repeat this extension process
until $R_{m+1}-T_m$ has been extended to a trivial constant rank projection $T_m^{'}$ defined on
$V_m \cap (V_{m+1} \cup \cdots \cup V_s)$ and satisfying
\[
T_m^{'}(x) \leq \left( \bigvee_{m < l \leq s} R_l(x) \right) - T_m(x), \ \forall x \in V_m \cap (V_{m+l} \cup \cdots \cup V_s).
\]
Note that the rank of $T_m^{'}$ is at least $d$.
Choose a trivial subprojection $T_m^{''}$ of $T_m^{'}$ with the property that
$\mathrm{rank}(T_m \oplus T_m^{''}) = \mathrm{rank}(R_{m+1}) - d$.  Apply Lemma
\ref{extendone} to extend $T_m^{''}$ to $V_m$ inside the complement of $T_m$,
and put
\[
R_m(x) = T_m(x) \oplus T_m^{''}(x), \ \forall x \in V_m.
\]
$R_m(x)$ so defined has the desired properties.

\end{proof}

\begin{lms}\label{trivbasis}
Let $X$ be a finite simplicial complex of dimension $d$, and let $U_1,\ldots,U_k$ be subcomplexes.
Suppose that for each $1 \leq i \leq k$, there exists a constant rank projection-valued function
$R_i \in \mathrm{M}_{\infty}(\mathrm{C}(U_i))$ satisfying:
\begin{enumerate}
\item[(i)] $R_i$ corresponds to a trivial vector bundle over $U_i$;
\item[(ii)] $\mathrm{rank}(R_1) \geq \lceil d/2 \rceil +1$;
\item[(iii)] $\mathrm{rank}(R_i) + \lceil d/2 \rceil +1 \leq \mathrm{rank}(R_{i+1})$, $1 \leq i < k$;
\item[(iv)] $R_i(x) \leq R_j(x)$ whenever $x \in U_i \cap U_j$ and $i \leq j$.
\end{enumerate}
Finally, suppose that $\mathrm{rank}(R_i) = n_i$, $n_1 < \cdots < n_k$, and
put $n_0=0$.  Then, 
there exist mutually orthogonal rank one projections 
$p_1,\ldots,p_{n_k} \in \mathrm{M}_{\infty}(\mathrm{C}(X))$,
each corresponding to a trivial vector bundle on $X$, such that
\[
R = \bigvee_{1 \leq i \leq k} R_i = \bigoplus_{i=1}^k \left(
\bigoplus_{j=n_{i-1}+1}^{n_i} \chi(U_i \cup \cdots \cup U_k) p_j \right).
\]
\end{lms}

\begin{proof}
We proceed by induction on $k$.  Consider the case $k=1$.  Since we
may view $R_1$ as upper semicontinuous projection-valued map from
$X$ into $\mathrm{M}_n(\mathbb{C})$ for
$n$ large, and since $(X,U_1)$ is a relative cell complex with cells
of dimension at most $d \ll n$, we may apply Lemma \ref{extendone}
to extend $R_1$ to a projection $\widetilde{R}_1$ on $X$ corresponding
to a trivial vector bundle.  Write $\widetilde{R}_1 = p_1 + \cdots + p_{n_1}$
for line bundles $p_1,\ldots,p_{n_1}$, each of which corresponds to a trivial
vector bundle.  Then,
\[
R_1 = \chi(U_1) \widetilde{R}_1 = \bigoplus_{j=1}^{n_1} \chi(U_1) p_j,
\]
as desired.

Now suppose that the lemma holds for $i < k$.  Conditions (i) and (iii)
in the hypotheses of the lemma together with Lemma \ref{extendone} 
allow us to extend $R_1$ to a trivial projection on $U_1 \cup U_2$ subordinate to $R_1 \vee R_2$.  
Iterating this process, we extend $R_1$ to a trivial projection $\widetilde{R}_1$
on $U_1 \cup \cdots \cup U_k$ which is subordinate to $R_1 \vee \cdots \vee R_k$.
Finally, apply Lemma \ref{extendone} once more to extend $\widetilde{R}_1$ to a trivial
projection on all of $X$.  Write $\widetilde{R}_1 = p_1 + \cdots + p_{n_1}$, where
each $p_j$ is projection corresponding to a one-dimensional trivial vector
bundle on $X$.  Then, as before,
\[
R_1 = \chi(U_1) \widetilde{R}_1 = \bigoplus_{j=1}^{n_1} \chi(U_1) p_j.
\] 

By \cite[Theorem 1.5, Chapter 8]{H}, we have that
$R_i - R_1$ is trivial for $1 \leq i \leq k$.  If $i > 1$, then $R_i-R_1$ has
rank greater than or equal to $\lceil d/2 \rceil +1$.  It follows that the projections
$(R_2-\widetilde{R}_1),(R_3-\widetilde{R}_1),\ldots,(R_k-\widetilde{R}_1)$ over the subcomplexes  
$U_2,U_3,\ldots,U_k$, respectively, satisfy the hypotheses of the lemma.  Moreover,
these projections may be viewed as maps from $X$ into the orthgonal complement of 
$\widetilde{R}_1$.  It follows that there exist rank one projections $p_{n_1+1},\ldots,p_{n_k} \in
\mathrm{M}_{\infty}(\mathrm{C}(X))$, each corresponding to a trivial vector bundle and
each orthogonal to $\widetilde{R}_1$, which satisfy
\[
\bigvee_{2 \leq i \leq k} (R_i-\widetilde{R}_1) = \bigoplus_{i=2}^k \left(
\bigoplus_{j=n_{i-1}+1}^{n_i} \chi(U_i \cup \cdots \cup U_k) p_j \right).
\]
We have
\[
R = \bigvee_{1 \leq i \leq k} R_i = R_1 \vee \left(\bigvee_{2 \leq i \leq k} (R_i-\widetilde{R}_1)\right),
\]
and the lemma follows.
\end{proof}

\begin{thms}\label{trivmajor}
Let $X$ be a finite simplicial complex, and let $a \in \mathrm{M}_{\infty}(\mathrm{C}(X))_+$.
Then, $a \precsim d$ for any $d \in \mathrm{M}_{\infty}(\mathrm{C}(X))_+$ which is trivial 
and satisfies 
\[
\mathrm{rank}(a)(x) + 4\mathrm{dim}(X) +3 \leq \mathrm{rank}(d)(x), \ \forall x \in X.
\]
\end{thms}

\begin{proof}
We may assume that $||a|| \leq 1$.  Let $n_1 < \cdots < n_k$ be the rank values taken by $a$, and
let $\epsilon>0$ be given. Form an approximant $f$ to $a$ satisfying the
conclusions of Theorem \ref{wellsupapprox}, where the sets $\overline{F_1},\ldots,\overline{F_k}$
corresponding to $f$ (cf. Definition \ref{trivial}) are subcomplexes of $X$
satisfying the conlusion of Lemma \ref{cellapprox}.  Let $p_1,\ldots,p_k$ be the
supporting projections for $f$, and notice that these satisfy the hypotheses of
Lemma \ref{trivcomp}.  The conclusion of Lemma \ref{trivcomp} then provides a family of 
constant rank projections $R_l \in \mathrm{M}_{\infty}(\mathrm{C}(V_l))$, $1 \leq l \leq s$,
(each $V_l$ is a union of consecutive $\overline{F_i}$s).  These, in addition to satisfying
the hypotheses of Lemma \ref{trivbasis}, have the properties that
\[
R(x) := \bigvee_{1 \leq l \leq s} R_l(x) \geq \bigvee_{1 \leq i \leq k} p_i(x) \geq f(x), \ \forall x \in X,
\]
and
\[
\mathrm{rank}(R)(x) - \mathrm{rank}\left(\bigvee_{1 \leq i \leq k} p_i(x) \right) \leq 4\mathrm{dim}(X) +3.
\]
Set $m_0=0$ and $m_l = \mathrm{rank}(R_l)$.
Let $q_1,\ldots,q_{n_k}$ be the family of mutually orthogonal rank one projections,
each corresponding to a trivial bundle, which are provided by the conclusion of Lemma
\ref{trivbasis}, i.e.,
\[
R = \bigoplus_{l=1}^s \left(
\bigoplus_{j=m_{l-1}+1}^{m_l} \chi(V_l \cup \cdots \cup V_s) q_j \right).
\]

By construction (see the proof of Theorem \ref{wellsupapprox}) we have that 
\[
H_{i,a} = \{ x \in X| \mathrm{rank}(a)(x) \leq n_i\} \subseteq (\overline{F_1} 
\cup \cdots \cup \overline{F_i})^{\circ}.
\]
Put
\[
F_{i,a} = \{ x \in X| \mathrm{rank}(a)(x) = n_i \}.
\]
Use Lemma \ref{cellapprox} to find, inductively (and upon refining the simplicial structure of $X$ if
necessary), subcomplexes $U_1,\ldots,U_k$ of $X$ satisfying:
\[
\overline{F_i} \cup \cdots \cup \overline{F_k} \subseteq (U_i \cup \cdots \cup U_k)^{\circ} 
\subseteq U_i \cup \cdots \cup U_k  \subseteq H_{i-1,a}^c, \ 1 \leq i \leq k
\]
The set $V_l$ is a union of consecutive $\overline{F_i}$s;  let $J_l$ be the set
of indices occurring among these $\overline{F_i}$s, and let $M_l$ be the largest
element of $J_l$.  Put 
\[
\widetilde{V}_l = \bigcup_{i \in J_l} U_i.
\]
For each $1 \leq l \leq s$, choose a continuous function $g_l:X \to [0,1]$ which
is identically one on $V_l \cup \cdots \cup V_s \subseteq H_{M_{l-1},a}^c$, identically
zero on $H_{M_{l-1},a}$, and nonzero off $H_{M_{l-1},a}$.  Then,
\[
\widetilde{R}(x) = \bigoplus_{l=1}^s \left(\bigoplus_{j=m_{l-1}+1}^{m_l} g_l(x) q_j \right)
\] 
is a trivial element of $\mathrm{M}_{\infty}(\mathrm{C}(X))_+$ such that $\widetilde{R}(x) \geq
R(x)$, $\forall x \in X$.

We claim that for all $x \in X$, we have 
\[
\mathrm{rank}(\widetilde{R})(x) - \mathrm{rank}(a)(x) \leq 3 \mathrm{dim}(X) +3.
\]
Indeed, 
\[
A_l := \{x \in X| \mathrm{rank}(\widetilde{R})(x) = m_l\} = V_l \cap (V_1 \cup \cdots V_{l-1})^c,
\]
and $g_j$ is nonzero on $A_l$ if and only if $j \leq l$.  We have chosen $g_j$
to satisfy
\[
g_j(x) \neq 0 \Longleftrightarrow x \in H_{M_{j-1},a}^c,
\] 
so $x \in A_l$ if and only if
\[
x \in H_{M_{j-l},a}^c \cap H_{M_l,a}.
\]
Thus, for $x \in A_l$, $\mathrm{rank}(a)(x) = n_i$ for some $i \in J_l$.
We have that
\[
\mathrm{rank}(R)(x) - \mathrm{rank}\left(\bigvee_{1 \leq i \leq k} p_i(x) \right) \leq 4 \mathrm{dim}(X) +3,
\]
whence $m_l - n_i \leq 4 \mathrm{dim}(X) + 3$, $\forall i \in J_l$. This proves the claim.

If $d$ is trivial and satisfies the hypotheses of the theorem, then
\[
\mathrm{rank}(d)(x) - \mathrm{rank}(\widetilde{R})(x) \geq 0, \ \forall x \in X.
\]
It follows from Proposition \ref{trivrankcomp} that 
\[
f \leq \widetilde{R} \precsim d.
\]
Thus, for every $\epsilon>0$, there exists $v \in \mathrm{M}_{\infty}(\mathrm{C}(X))$
such that 
\[
||vdv^*-f||<\epsilon.  
\]
It follows that
\[
||vdv^* - a|| \leq \epsilon + ||vdv^*- f|| < 2 \epsilon;
\]
$\epsilon$ was arbitrary, and the theorem follows.
\end{proof}

\subsection{Trivial minorants}
\begin{thms}\label{trivminor}
Let $X$ be a compact metric space, and $a \in \mathrm{M}_{\infty}(\mathrm{C}(X))_+$.
Then, $d \precsim a$ whenever $d$ is trivial and satisfies
\[
\mathrm{rank}(d)(x) \leq \mathrm{max}\{\mathrm{rank}(a)(x)-\mathrm{dim}(X)-1,0\}, \ \forall x \in X.
\]
\end{thms}

\begin{proof}  The hypotheses of the theorem imply that $d = 0$ if $\mathrm{dim}(X) = \infty$,
in which case the theorem holds.  

Suppose that $\mathrm{dim}(X) < \infty$.
We proceed by induction on $k$.
Suppose that $a$ takes rank values $n_1 < \cdots < n_k$, and, with $n_0 = 0$, put 
\[
G_i = \{x \in X|\mathrm{rank}(a)(x) > n_{i-1} \}.
\] 

Suppose that $k=1$.  If $n_1 \leq \mathrm{dim}(X)+1$, then there is nothing to prove, so
assume that $n_1 > \mathrm{dim}(X)+1$.  By \cite[Theorem 1.5, Chapter 8]{H}, $a$ (which,
since it only takes one rank value, is Cuntz equivalent to a projection) is Cuntz
equivalent to $\theta_{n_1-\mathrm{dim}(X)+1} \oplus p$ for some projection $p \in
\mathrm{M}_{\infty}(\mathrm{C}(X))$.  The projection $\theta_{n_1-\mathrm{dim}(X)+1}$
is trivial, and so dominates any trivial $d$ satisfying the hypotheses of the theorem
by Proposition \ref{trivrankcomp}.

Now suppose that we have proved the theorem when $a$ takes $i < k$ rank values.
We treat two cases:  $n_1 \leq \mathrm{dim}(X)+1$ and $n_1 > \mathrm{dim}(X)+1$.  

Suppose first that $n_1 \leq \mathrm{dim}(X)+1$.  Then, any trivial $d$ satisfying
the hypotheses of the theorem necessarily satisfies 
\[
\mathrm{rank}(d)(x) = 0, \ \forall x \in G_2^c.
\]
Put $\epsilon_n = 1/2^n$.  Since $||d(x)||$ is uniformly continuous on $\overline{G_2}$, there exists
a sequence $(\delta_n)$ of nonnegative reals such that $||d(x)-d(y)||< \epsilon_n/3$
whenever $\mathrm{dist}(x,y) < \delta_n$.  For each $n \in \mathbb{N}$, let
$V_n \subseteq \overline{G_2}$ be the closed set which is the complement of the
set of points in $\overline{G_2}$ whose distance from $\overline{G_2} \backslash G_2$
is strictly less than $\delta_n$, and let $U_n$ be the closed set consisting of
those points in $\overline{G_2}$ whose distance from $\overline{G_2} \backslash G_2$
is less than or equal to $\delta_n/2$.  
Choose a function $f_n:\overline{G_2} \to [0,1]$
which is identically zero on $U_n$ and identically one
on $V_n$.  Notice that $||d(x)|| < \epsilon$, $\forall x \in \overline{G_2} \backslash V_n$.  
Upon restriction to $\overline{U_n^c}$, $a$ takes at most $k-1$ rank values. Since $a|_{\overline{U_n^c}}$ and $d|_{\overline{U_n^c}}$
satisfy the hypothese of the theorem {\it a fortiori}, there is an element $w_n \in \mathrm{M}_{\infty}(\mathrm{C}(\overline{U_n^c}))$
such that 
\[
||w_n(a|_{\overline{U_n^c}})w_n^* - d|_{\overline{U_n^c}}|| < \epsilon_n.
\]
Put $g_n = f_n \cdot w_n$, and note that since $f_n$ is zero off $\overline{U_n^c}$, we may
view $g_n$ as an element of $\mathrm{M}_{\infty}(\mathrm{C}(\overline{G_2}))$.  Now, 
\[
||g_n(x)a(x)g_n^*(x) - d(x) || < \epsilon_n, \ \forall x \in \overline{V_n},
\]
and for every $x \in V_n^c$
\begin{eqnarray*}
||g_n(x)a(x)g_n^*(x) - d(x) || & \leq & ||f_n(x)w_n(x)a(x)w_n^*(x)f_n(x)|| + ||d(x)|| \\
& \leq & ||f_n(x)^2(w_n(x)a(x)w_n^*(x) - d(x))|| + 2||d(x)|| \\
& \leq & 3||d(x)|| < \epsilon_n
\end{eqnarray*}
(the second-to-last inequality uses that $f_n(x) = 0$, $\forall x \in U_n$.) 
Thus,
\[
g_n a g_n^* \stackrel{n \to \infty}{\longrightarrow} d,
\]
and $d \precsim a$, as desired.

Now suppose that $n_1 > \mathrm{dim}(X)+1$.  We will reduce to the case $n_1 \leq \mathrm{dim}(X)+1$,
proving the theorem.  By \cite[Proposition 3.2]{DNNP}, there is a projection $p$ on $X$ such
that 
\[
\mathrm{rank}(p) \geq \lceil \mathrm{dim}(X)/2 \rceil +(n_1 - \mathrm{dim}(X)-1)
\]
and $p \precsim a$.  An application of \cite[Theorem 1.2, Chapter 8]{H} yields a trivial
subprojection $\theta_m$ of $p$, where $m = n_1 - \mathrm{dim}(X) - 1$.  By
\cite[Proposition 2.2]{PT}, there is a positive element $b \in \mathrm{M}_{\infty}(\mathrm{C}(X))$
such that $a \sim b \oplus \theta_m$.  Let $p_1,\ldots,p_m$ be the first $m$ trivial
rank one projections supporting $d$ (cf. Definition \ref{trivial}).  Then,
$p_1 \oplus \cdots \oplus p_m \sim \theta_m$ (where $\sim$ is in fact true Murray-von Neumann
equivalence, and 
\[
d^{'} := (1-(p_1 \oplus \cdots \oplus p_m))d(1-(p_1 \oplus \cdots \oplus p_m))
\]
and
\[
d^{''} := (p_1 \oplus \cdots \oplus p_m)d(p_1 \oplus \cdots \oplus p_m)
\]
are trivial.  Moreover, $d^{'}$ and $b$ (substituted for $d$ and $a$, respectively)
satisfy the hypotheses of the theorem, $b$ takes $k$ rank values, and the lowest
rank value taken by $b$ is less than or equal to $\mathrm{dim}(X)+1$.  We may thus
apply our proof above to conclude that $d^{'} \precsim b$.  Since $d^{''} \precsim \theta_m$
we have 
\[
d = d^{'} \oplus d^{''} \precsim b \oplus \theta_m \sim a,
\]
as desired.
\end{proof}

\subsection{The main theorem}

\begin{thms}\label{main}
Let $X$ be a compact metric space of covering dimension $d \in \mathbb{N}$.  
Let $a,b \in \mathrm{M}_n(\mathrm{C}(X))$ be positive, and suppose that 
\[
\mathrm{rank}(a)(x) + 9 d \leq \mathrm{rank}(b)(x), \ \forall x \in X.
\]
Then, $a \precsim b$. 
\end{thms}

\begin{proof}
Combining Theorems \ref{trivmajor} and \ref{trivminor} yields the theorem
with $9d$ replaced by $5d+4$.  But if $d=0$, then
\[
\mathrm{rank}(a)(x) \leq \mathrm{rank}(b)(x), \ \forall x \in X
\]
implies $a \precsim b$ by \cite[Theorem 3.3]{Pe1}.  If $d \geq 1$, then
$5d+4 \leq 9d$.
\end{proof}

Theorem \ref{main} applies equally well to locally compact second countable
Hausdorff spaces whose one point compactifications have finite covering dimension.
We fully expect that it will generalise to recursive subhomogeneous $C^*$-algebras.



\section{Applications to AH algebras}\label{ah}

\subsection{The dimension-rank ratio vs. the radius of comparison}

Recall the terminology concerning AH algebras from Section 1.
A unital AH algebra $A$ has \emph{flat dimension growth} (\cite[Definition 1.2]{To3}) if
it admits a decomposition for which
\begin{equation}\label{moddimgrowth}
\limsup_{i \to \infty} \ \mathrm{max}_{1 \leq l \leq m_i} \left\{ 
\frac{\mathrm{dim}(X_{i,l})}{\mathrm{rank}(p_{i,l})} \right\} < \infty
\end{equation}
To study unital AH algebras with flat dimension growth was suggested
by Blackadar (\cite{Bl2}) in 1991, but there were no non-trivial examples of
such algebras --- algebras with flat dimension growth but {\it not} slow dimension
growth --- until the pioneering work of Villadsen in 1997 (\cite{V1}).
We initiated the study of such algebras in earnest in \cite{To3}.  Our key
tool was the {\it dimension-rank ratio} of a unital AH algebra $A$ (write $\mathrm{drr}(A)$), an
isomorphism invariant which is defined to be the infimum of the set of 
nonnegative reals $c$ such that $A$ has a decomposition satisfying
\begin{equation}\label{dimrank}
\limsup_{i \to \infty} \ \mathrm{max}_{1 \leq l \leq m_i} \left\{ 
\frac{\mathrm{dim}(X_{i,l})}{\mathrm{rank}(p_{i,l})} \right\} \leq c.
\end{equation}
This invariant may be thought of as a measure of the ratio of the topological dimension
of $A$ to its matricial size, despite the fact that both quantities may be
infinite.  Its effectiveness in studying AH algebras of flat dimension growth
leads one naturally to consider whether it has an analogue defined for all
unital and stably finite $C^*$-algebras.  In \cite{To3} we introduced the following 
candidate:

\begin{dfs}[Definition 6.1, \cite{To3}]\label{rc}
Say that $A$ has $r$-comparison if whenever one has positive elements
$a,b \in \mathrm{M}_{\infty}(A)$ such that
\[
s(\langle a \rangle) + r < s(\langle b \rangle), \ \forall s \in \mathrm{LDF}(A),
\]
then $\langle a \rangle \leq \langle b \rangle$ in $W(A)$.
Define the radius of comparison
of $A$, denoted $\mathrm{rc}(A)$, to be 
\[
\mathrm{inf} \{r \in \mathbb{R}^+ | \ (W(A), \langle 1_A \rangle) \ \mathrm{has} \ r-\mathrm{comparison} \ \}
\]
if it exists, and $\infty$ otherwise.
\end{dfs}
\noindent 

Theorem \ref{main} confirms the radius of comparison as the proper abstraction of
the dimension-rank ratio for semi-homogeneous algebras.   

\begin{thms}\label{semihomrc}
There exist constants $K_1, K_2 > 0$ such that for
any semi-homogeneous $C^*$-algebra
\[
A = \bigoplus_{i=1}^n p_i(\mathrm{C}(X_i) \otimes \mathcal{K})p_i
\]
with each $X_i$ a connected finite-dimensional CW-complex one has the following inequalities:
\[
\mathrm{drr}(A) \leq K_1 \mathrm{rc}(A) 
\]
and
\[
\mathrm{rc}(A) \leq K_2 \mathrm{drr}(A).
\]
\end{thms}

\begin{proof}  Suppose first that $n=1$.
The first inequality is \cite[Theorem 6.6]{To3}.  For the second inequality, use \cite[Theorem 2.2]{To3}
to conclude that $\mathrm{drr}(A) = \mathrm{dim}(X)/\mathrm{rank}(p)$;  the inequality now follows
from the definition of the radius of comparison and Theorem \ref{main}.  To prove the theorem
for general $n$ use the following facts:
\[
\mathrm{drr}(A \oplus B) = \mathrm{max} \{\mathrm{drr}(A),\mathrm{drr}(B)\}
\]
and
\[
\mathrm{rc}(A \oplus B) = \mathrm{max} \{\mathrm{rc}(A),\mathrm{rc}(B)\}
\]
for any unital AH algebras $A$ and $B$ (\cite[Proposition 2.2 (ii)]{To3} and
\cite[Proposition 6.2 (ii)]{To3}, resp.). 
\end{proof}

We will address the relationship between the dimension-rank ratio and the radius of comparison
for general AH algebras in a separate paper.

\subsection{Strict comparison of positive elements}

The next lemma is due to M. R{\o}rdam.  We are grateful for his permission to use it here.

\begin{lms}[R{\o}rdam, \cite{R7}]\label{posapprox1}
Let $A$ be a $C^*$-algebra and $\{A_i\}_{i \in I}$ a collection of $C^*$-subalgebras whose
union is dense.   Then, for every $a \in \mathrm{M}_{\infty}(A)_+$ and $\epsilon > 0$
there exists $i \in I$ and $\tilde{a} \in \mathrm{M}_{\infty}(A_i)$ such that
\[
(a-\epsilon)_+ \precsim \tilde{a} \precsim (a-\epsilon/2)_+ \precsim a 
\]
in $\mathrm{M}_{\infty}(A)$
\end{lms}

\begin{proof}
First find a positive element $b$ in some $\mathrm{M}_{\infty}(A_i)$, $i \in I$,
such that 
\[
||b-(a-\epsilon/2)_+||< \epsilon/4.  
\]
Put $\tilde{a} := (b-\epsilon/4)_+$.  
The conclusion follows from Proposition \ref{basics} and the estimate $||a-\tilde{a}|| < \epsilon$.
\end{proof}

\begin{lms}\label{posapprox2}
Let $A$ be the limit of an inductive system $(A_i,\phi_i)_{i \in \mathbb{N}}$
of $C^*$-algebras, where $\phi_i$ is injective for each $i \in \mathbb{N}$.  
Let $a,b \in \mathrm{M}_{\infty}(A_i)$ be positive elements
such that $\phi_{i \infty}(a) \precsim \phi_{i \infty}(b)$ in $\mathrm{M}_{\infty}(A)$.
Then, for every $\epsilon>0$ there is a $j>i$ such that 
\[
(\phi_{ij}(a)-\epsilon)_+ \precsim \phi_{ij}(b)
\]
inside $\mathrm{M}_{\infty}(A_j)$.
\end{lms}

\begin{proof}
By working in a matrix algebra over $A$, we may assume that $a,b \in A$.
Since the $\phi_i$ are injective, we simply identify $a$ and $b$ with their
forward images in $A_j$, $j \geq i$, and in $A$ itself.  We have 
$a \precsim b$ in $\mathrm{M}_{\infty}(A)$, so there is a sequence
$(v_n)$ in $\mathrm{M}_{\infty}(A)$ such that 
\[
v_nbv_n^* \stackrel{n \to \infty}{\longrightarrow} a.
\]
This sequence may be chosen to lie in the dense local $C^*$-algebra
$\cup_{i=1}^{\infty} A_i$.  Indeed, for any $w_n \in A$ we have
\begin{eqnarray*}
||w_n b w_n* - v_n b v_n^*|| & = & ||(w_n-v_n+v_n)b(w_n-v_n+v_n)^* - v_nbv_n^*||\\
& = & ||(w_n-v_n)b(w_n-v_n)^* + (w_n-v_n)bv_n* + v_nb(w_n-v_n)^* || \\
& \leq & ||(w_n-v_n)||(||b(w_n-v_n)^*|| + ||bv_n^*|| + ||bv_n||),
\end{eqnarray*}
so choosing $w_n \in \cup_{i=1}^{\infty} A_i$ sufficiently close to $v_n$ yields
\[
w_nbw_n^* \stackrel{n \to \infty}{\longrightarrow} a.
\]
Let $\epsilon>0$ be given.  Find $i,n \in \mathbb{N}$ such that
\[
||w_nbw_n^* - a|| < \epsilon
\]
and $a,b,w_n \in \mathrm{M}_{\infty}A_i$.  It then follows from Proposition \ref{basics} that
$(a-\epsilon)_+ \precsim b$ inside $\mathrm{M}_{\infty}(A_i)$, as desired.
\end{proof}

\begin{thms}\label{rczerolim}
Let $(A_i,\phi_i)$ be an inductive sequence of unital, exact, and stably finite $C^*$-algebras
with simple limit $A$.  Suppose further that each $\phi_i$ is injective and that 
\[
\liminf_{i \to \infty} \mathrm{rc}(A_i) = 0.
\]
Then, $\mathrm{rc}(A)=0$.  In particular, $A$ has strict comparison of positive elements.
\end{thms}

\begin{proof}
We will prove that $W(A)$ is almost unperforated, i.e., that whenever one
has $\langle a \rangle, \langle b \rangle \in W(A)$ such that $(n+1)\langle a \rangle
 \leq n \langle b \rangle$, then $\langle a \rangle \leq \langle b \rangle$. 
It then follows from \cite[Corollary 4.6]{R4} that $A$ has strict comparison 
of positive elements. 

Let $a,b \in \mathrm{M}_{\infty}(A)$ be positive, and suppose that $(n+1)a \precsim nb$
for some $n \in \mathbb{N}$.  Let $\epsilon>0$ be given, and find $\delta>0$ such
that 
\[
((n+1)a-\epsilon/2)_+ = (n+1)(a-\epsilon/2)_+ \precsim n(b-\delta)_+ = (nb-\delta)_+.
\]
Use Lemma \ref{posapprox1} to find some $i \in \mathbb{N}$ and $\tilde{a},\tilde{b} \in
\mathrm{M}_{\infty}(A_i)_+$ such that: 
\begin{equation}\label{aest}
(a-3\epsilon/4)_+ \precsim \tilde{a} \precsim (a-\epsilon/2)_+; \ \ (b-\delta)_+ \precsim \tilde{b} \precsim {b}.
\end{equation}
It follows that
\[
(n+1)(a-\epsilon)_+ \precsim (n+1)(\tilde{a} - \epsilon/4)_+ \precsim n\tilde{b} \precsim nb.
\]
By Lemma \ref{posapprox2} we may, by increasing $i$ if necessary, assume that 
\begin{equation}\label{ineq3}
(n+1)(\tilde{a}-\epsilon/4)_+ \precsim n\tilde{b}
\end{equation}
inside $\mathrm{M}_{\infty}(A_i)$.  Since $A$ is simple, we may assume that the images
of both $(\tilde{a}-\epsilon/4)_+$ and $\tilde{b}$ under any $\tau \in \mathrm{T}(A_i)$
are non-zero.  It follows that 
\[
s((\tilde{a}-\epsilon/4)_+) \neq 0, \ s(\tilde{b}) \neq 0, \ \forall s \in \mathrm{LDF}(A_i).
\]
Equation (\ref{ineq3}) shows that for any $\tau \in \mathrm{T}(A_i)$,
\[
s_{\tau}(\tilde{b}) - s_{\tau}((\tilde{a}-\epsilon/4)_+) \geq (1/n) s_{\tau}(\tilde{b}).
\]
The map $\tau \mapsto s_{\tau}(\tilde{b})$ is strictly positive and lower semicontinuous
on the compact space $\mathrm{T}(A_i)$.  It therefore achieves a minimum value $c>0$, and
\[
s_{\tau}((\tilde{a}-\epsilon/4)_+) + c/2 < s_{\tau}(\tilde{b}), \ \forall \tau \in \mathrm{T}(A_i).
\]
Increasing $i$ if necessary, we may assume that $\mathrm{rc}(A_i) < c/2$, whence
\begin{equation}\label{best}
(\tilde{a}-\epsilon/4)_+ \precsim \tilde{b} \precsim b
\end{equation}
by the definition of $\mathrm{rc}(\bullet)$.  From (\ref{aest}) we have the inequality
\[
(a-3\epsilon/4)_+ \precsim \tilde{a},
\]
whence
\begin{equation}\label{aest2}
(a-\epsilon)_+ \precsim (\tilde{a} -\epsilon/4)_+
\end{equation}
by the functional calculus.  Combining (\ref{best}) and (\ref{aest2}) we have
\[
(a-\epsilon)_+ \precsim b.
\]
Since $\epsilon$ was arbitrary, the theorem follows from Proposition \ref{basics}.
\end{proof}

\begin{cors}\label{ahaup}
Let $A$ be a simple unital AH algebra with slow dimension 
growth.  Then, $W(A)$ is almost unperforated.  In particular, $A$
has strict comparison of positive elements.
\end{cors}

\begin{proof}
As in equation (\ref{decomp}) we have
$A= \lim_{i \to \infty}(A_i,\phi_i)$, where
\[
A_i = \bigoplus_{l=1}^{n_i} p_{i,l}(\mathrm{C}(X_{i,l}) \otimes \mathcal{K})p_{i,l}.
\]
$A$ has slow dimension growth, so we may assume that the $\phi_i$ are injective (this is the main result of
\cite{EGL3}) and that 
\[
\mathrm{max}_{1 \leq l \leq n_i} \ \left\{ \frac{\mathrm{dim}(X_{i,1})}{\mathrm{rank}(p_{i,1})}, \ldots,
\frac{\mathrm{dim}(X_{i,n_i})}{\mathrm{rank}(p_{i,n_i})} \right\} \stackrel{i \to \infty}{\longrightarrow} 0.
\]
This last condition, by \cite[Theorem 2.3]{To3}, implies that $\mathrm{drr}(A_i) \to 0$ as $i \to \infty$.
Proposition \ref{semihomrc} then shows that $\mathrm{rc}(A_i) \to 0$ as $i \to \infty$, and we
have collected the hypotheses of Theorem \ref{rczerolim}.
\end{proof}


Combining Corollary \ref{ahaup} with earlier work of R{\o}rdam we recover the characterisation
of real rank zero for simple unital AH algebras of slow dimension growth obtained in \cite{BDR}.

\begin{cors}[Blackadar, D\u{a}d\u{a}rlat, and R{\o}rdam, \cite{BDR}]
Let $A$ be a simple unital AH algebra with slow dimension growth.  Then, the following are equivalent:
\begin{enumerate}
\item[(i)] $A$ has real rank zero;
\item[(ii)] the projections in $A$ separate tracial states;
\item[(iii)] the image of $\mathrm{K}_0(A)$ is uniformly dense in the space of continuous affine
functions on the tracial state space $\mathrm{T}(A)$.
\end{enumerate}
\end{cors}

\begin{proof}
$A$ has stable rank one by the results of \cite{BDR}, and a weakly unperforated Cuntz semigroup
by Corollary \ref{ahaup}.  Apply \cite[Theorem 7.2]{R4}.
\end{proof}

Finally, we have that several theorems on $\mathcal{Z}$-stable $C^*$-algebras recently proved
by Brown, Perera, and the author can be extended to AH algebras with slow dimension growth.

\begin{cors}\label{conjectures}
Let $A$ be an element of $\mathcal{SDG}$.  The following statements hold:
\begin{enumerate}
\item[(i)] $W(A)$ is recovered functorially from the Elliott invariant of $A$;
\item[(ii)] $\mathrm{LDF}(A)$ is weak-$*$ dense in $\mathrm{DF}(A)$;
\item[(iii)] $\mathrm{DF}(A)$ is a Choquet simplex.
\end{enumerate}
\end{cors}

\begin{proof}
Since $A$ has strict comparison of positive elements, we may appeal to Theorem 5.4 and
Corollary 6.9 of \cite{BPT}.
\end{proof}

We refer the reader to Section 4 of \cite{PT}
for an explicit description of the functor which reconstructs $W(A)$ from the Elliott 
invariant of $A$ in (i) above.

\vspace{5mm}

\noindent
\emph{Andrew S. Toms \newline
Department of Mathematics and Statistics \newline
\hspace*{2mm} York University \newline
\hspace*{2mm} 4700 Keele St. \newline
\hspace*{2mm} Toronto, Ontario, M3J 1P3 \newline
Canada \newline}

\noindent
atoms@mathstat.yorku.ca

\end{document}